%% file: optsubspace.tex
\documentclass[review]{elsarticle}

\usepackage{lineno,hyperref}
\modulolinenumbers[5]

\journal{CMAME}

\usepackage{amsmath,amssymb}
\usepackage{graphicx} 
\usepackage{algorithm}
\usepackage{subfig} 
\usepackage{bm}

\input{pauls-macros}

\newcommand{\mCref}{\mC_{\text{ref}}}
\newcommand{\mWref}{\mW_{\text{ref}}}
\newcommand{\Lambdaref}{\Lambda_{\text{ref}}}

\newtheorem{theorem}{Theorem}
\newtheorem{definition}{Definition}
\newtheorem{assumption}{Assumption}

\begin{document}

\begin{frontmatter}

\title{A near-stationary subspace for ridge approximation}





\author[CSM]{Paul G.~Constantine} 
\ead{paul.constantine@mines.edu}

\author[TU]{Armin Eftekhari} 

\author[CSM]{Jeffrey Hokanson}

\author[UT]{Rachel A.~Ward}
\address[CSM]{Department of Applied Mathematics and Statistics, Colorado School of Mines, Golden, CO 80211}
\address[TU]{Alan Turing Institute, British Library, 96 Euston Road, London NW1 2DB}
\address[UT]{Department of Mathematics and Institute for Computational Engineering and Sciences, University of Texas at Austin, Austin, TX 78712}

\begin{abstract}
Response surfaces are common surrogates for expensive computer simulations in engineering analysis. However, the cost of fitting an accurate response surface increases exponentially as the number of model inputs increases, which leaves response surface construction intractable for high-dimensional, nonlinear models. We describe \emph{ridge approximation} for fitting response surfaces in several variables. A ridge function is constant along several directions in its domain, so fitting occurs on the coordinates of a low-dimensional subspace of the input space. We review essential theory for ridge approximation---e.g., the best mean-squared approximation and an optimal low-dimensional subspace---and we prove that the gradient-based \emph{active subspace} is near-stationary for the least-squares problem that defines an optimal subspace. Motivated by the theory, we propose a computational heuristic that uses an estimated active subspace as an initial guess for a ridge approximation fitting problem. We show a simple example where the heuristic fails, which reveals a type of function for which the proposed approach is inappropriate. We then propose a simple alternating heuristic for fitting a ridge function, and we demonstrate the effectiveness of the active subspace initial guess applied to an airfoil model of drag as a function of its 18 shape parameters. 
\end{abstract}


\end{frontmatter}

%\linenumbers

\input{sec0-intro}

\input{sec1-opt}

\input{sec2-subspaces}
\input{sec3-examples}
\input{sec4-conclusion}

\appendix

\input{sec-appendix3}

\section*{Acknowledgment}

\noindent The first author’s work is supported by the U.S. Department of Energy Office of Science, Office of Advanced Scientific Computing Research, Applied Mathematics program under Award Number DE-SC-0011077. The first and third authors' work is supported by Department of Defense, Defense Advanced Research Project Agency’s program Enabling Quantification of Uncertainty in Physical Systems. The second author's work is supported by the Alan Turing Institute under the EPSRC grant EP/N510129/1. The third author's work is partially funded by NSF CAREER grant \#1255631.  

\section*{References}

\bibliography{optsubspace}

\end{document}

%% file: pauls-macros.tex
\graphicspath{{./}{./figs/}}

  \def\clap#1{\hbox to 0pt{\hss#1\hss}}

\providecommand{\mat}[1]{\bm{#1}}%
\renewcommand{\vec}[1]{\mathbf{#1}}

\providecommand{\mC}{\ensuremath{\mat{C}}}

\providecommand{\mI}{\ensuremath{\mat{I}}}

\providecommand{\mQ}{\ensuremath{\mat{Q}}}

\providecommand{\mT}{\ensuremath{\mat{T}}}
\providecommand{\mU}{\ensuremath{\mat{U}}}
\providecommand{\mV}{\ensuremath{\mat{V}}}
\providecommand{\mW}{\ensuremath{\mat{W}}}

\providecommand{\ve}{\ensuremath{\vec{e}}}

\providecommand{\vg}{\ensuremath{\vec{g}}}

\providecommand{\vs}{\ensuremath{\vec{s}}}

\providecommand{\vu}{\ensuremath{\vec{u}}}
\providecommand{\vv}{\ensuremath{\vec{v}}}
\providecommand{\vw}{\ensuremath{\vec{w}}}
\providecommand{\vx}{\ensuremath{\vec{x}}}
\providecommand{\vy}{\ensuremath{\vec{y}}}
\providecommand{\vz}{\ensuremath{\vec{z}}}

\newcommand{\bmat}[1]{\begin{bmatrix}#1\end{bmatrix}}

\newcommand{\UyVz}{\mU\vy + \mV\vz}
\renewcommand{\span}{\operatorname{span}}
\newcommand{\pzgy}{\pi(\vz|\vy)}

\newcommand{\bnabla}{\bar{\nabla}}

\newcommand{\ppV}{\frac{\partial}{\partial \mV}}

\newcommand{\ppvij}{\frac{\partial}{\partial v_{ij}}}

%% file: sec0-intro.tex
\section{Introduction}
\label{sec:intro}

\noindent Engineering computations often employ cheap response surfaces that mimic the input/output relationship between an expensive computer model's parameters and its predictions. The essential idea is to use a few expensive model runs at particular parameter values (i.e., a \emph{design of experiments}~\cite{dace2003}) to fit or train a response surface, where the surface may be a polynomial, spline, or radial basis approximation~\cite{Jones2001,Barton1992}. The same scenario motivates statistical tools for design and analysis of computer experiments~\cite{Sacks1989,Kennedy2000,dace2003}, which use Gaussian processes to model uncertainty in the surrogate's predictions.

The cost of constructing an accurate response surface increases exponentially as the dimension of the input space increases; in approximation theory, this is the tractability problem~\cite{Wozniakowski1994,Traub1998}, though it is colloquially referred to as the \emph{curse of dimensionality}~\cite[Section 5.16]{bellman61}. Several techniques attempt to alleviate this curse---each with advantages and drawbacks for certain classes of problems; see~\cite{Shan2009} for an extensive survey. One idea is to identify unimportant input variables with global sensitivity metrics~\cite{saltelli2008global} and fix them at nominal values, which effectively reduces the dimension for response surface construction; Sobol' et al.~studied the effects of such coordinate-based dimension reduction on the approximation~\cite{Sobol2007}. A generalization of coordinate-based dimension reduction is to identify unimportant directions---not necessarily coordinate aligned. If the scientist can identify a few important linear combinations of inputs, then she may fit a response surface of only those linear combinations, which allows a higher degree of accuracy along important directions in the input space. 

A \emph{ridge function}~\cite{pinkus2015} is a function of a few linear combinations of inputs that takes the form $g(\mU^T\vx)$, where $\vx\in\mathbb{R}^m$, $\mU\in\mathbb{R}^{m\times n}$ with $n<m$, and $g:\mathbb{R}^n\rightarrow\mathbb{R}$. The term \emph{ridge function} is more commonly used when $\mU$ is a single vector ($n=1$). Pinkus calls our definition a \emph{generalized ridge function}~\cite[Chapter 1]{pinkus2015}, though Keiper uses the qualifier \emph{generalized} for a model where $\mU$ depends on $\vx$~\cite{Keiper2015}. A ridge function is constant along directions in its domain that are orthogonal to $\mU$'s columns. To see this, let $\vv\in\mathbb{R}^m$ be orthogonal to $\mU$'s columns; then
\begin{equation}
g(\mU^T(\vx + \vv)) = g(\mU^T\vx + \underbrace{\mU^T\vv}_{=\;0}) = g(\mU^T\vx).
\end{equation}
If $\mU$ is known, then one need only construct $g$, which is a function of $n<m$ variables. Thus, constructing $g$ may require exponentially fewer model evaluations than constructing a comparably accurate response surface on all $m$ variables. 

Let $f:\mathbb{R}^m\rightarrow\mathbb{R}$ represent the simulation model's input/output map to approximate, and let its domain be equipped with a probability function $\rho:\mathbb{R}^m\rightarrow\mathbb{R}_+$. The function $\rho$ may model uncertainty in the simulation's input parameters, which is a common modeling choice in \emph{uncertainty quantification}~\cite{Smith2013,Sullivan2015}.
The ridge approximation problem may be stated as: given $f$ and $\rho$, find $g$ and $\mU$ that minimize the approximation error. After a brief survey of related concepts, we define a specific ridge approximation problem in Section \ref{sec:optU}. We then study a particular $\mU$ derived from $f$'s gradient known as the \emph{active subspace}~\cite{asmbook}. We show that, under certain conditions, the active subspace is nearly stationary---i.e., that the gradient of the objective function defining the approximation problem is bounded; see Section \ref{sec:stationary}. This result motivates a heuristic for the initial subspace when fitting a ridge approximation given pairs $\{(\vx_i,f(\vx_i))\}$. In Section \ref{sec:examples}, we show a simple bivariate example that exposes the limitations of the heuristic. We then study an 18-dimensional example from an airfoil shape optimization problem where the heuristic succeeds; in particular, we demonstrate a numerical procedure for estimating the active subspace using samples of the gradient, and we show how the estimated active subspace is a superior starting point for a numerical optimization heuristic for fitting the ridge approximation. 

\subsection{Related concepts}

\noindent There are many concepts across subfields that relate to ridge approximation. In what follows, we briefly review three of these subfields with citations that point interested readers to representative works. 

\subsubsection{Projection pursuit regression}

\noindent In the context of statistical regression, Friedman and Stuetzle~\cite{Friedman1981} proposed \emph{projection pursuit regression} with a ridge function model:
\begin{equation}
\label{eq:ppr}
y_i \;=\; \sum_{k=1}^r g_k(\vu_k^T\vx_i) + \varepsilon_i,
\end{equation}
where $\vx_i$'s are samples of the predictors, $y_i$'s are the associated responses, and $\varepsilon_i$'s model random noise---all standard elements of statistical regression~\cite{Weisberg2005}. The $g_k$'s are smooth univariate functions (e.g., splines), and the $\vu_k$'s are the directions of the ridge approximation. To fit the projection pursuit regression model, one minimizes the mean-squared error over the directions $\{\vu_k\}$ and the parameters of $\{g_k\}$. Motivated by the projection pursuit regression model, Diaconis and Shahshahani~\cite{Diaconis1984} studied the approximation properties of nonlinear functions ($g_k$ in \eqref{eq:ppr}) of linear combinations of the variables ($\vu_k^T\vx$ in \eqref{eq:ppr}). Huber~\cite{Huber1985} surveyed a wide class of projection pursuit approaches across an array of multivariate problems; by his terminology, \emph{ridge approximation} could be called \emph{projection pursuit approximation}. Chapter 11 of Hastie, Tibshirani, and Friedman~\cite{hastie2009elements} links projection pursuit regression to neural networks, which uses ridge functions with particular choices for the $g_k$'s (e.g., the sigmoid function). Although algorithm implementations may be similar, the statistical regression context is different from the approximation context, since there is no inherent randomness in the approximation problem. 

\subsubsection{Gaussian processes with low-rank correlation models}

\noindent In Gaussian process regression~\cite{gpml2006}, the conditional mean of the Gaussian process model given data (e.g., $\{y_i\}$ as in \eqref{eq:ppr}) is the model's prediction. This conditional mean is a linear combination of radial basis functions with centers at a set of points $\{\vx_i\}$, where the form of the basis function is related to the Gaussian process' assumed correlation. Vivarelli and Williams~\cite{vivarelli99} proposed a correlation model of the form
\begin{equation}
C(\vx,\vx') \;\propto\;
\exp\left[
-\frac{1}{2}(\vx-\vx')^T\mU\mU^T(\vx-\vx')
\right],
\end{equation}
where $\mU$ is a tall matrix. In effect, the resulting conditional mean is a function of linear combinations of the predictors, $\mU^T\vx$---i.e., a ridge function. A maximum likelihood estimate of $\mU$ is the minimizer of an optimization similar to the one we define for ridge approximation; see Section \ref{sec:optU}. Bilionis et al.~\cite{bilionis2016gaussian}, use a similar approach from a Bayesian perspective in the context of uncertainty quantification, where the subspace defined by $\mU$ enables powerful dimension reduction. 

\subsubsection{Ridge function recovery}

\noindent Recent work in constructive approximation seeks to recover the parameters of a ridge function from point queries~\cite{Fornasier2012,cohen2012,Tyagi2014}. In other words, assume $f(\vx)=g(\mU^T\vx)$ is a ridge function; using pairs $\{\vx_i,f(\vx_i)\}$, one wishes to recover the components of $\mU$. Algorithms for determining $\mU$ (e.g., Algorithm 2 in~\cite{Fornasier2012}) are quite different than optimizing a ridge approximation over $\mU$. However, the recovery problem is similar in spirit to the ridge approximation problem.

%% file: sec1-opt.tex
\section{Optimal ridge approximation}
\label{sec:optU}

\noindent Consider a function $f:\mathbb{R}^m\rightarrow\mathbb{R}$ that is square-integrable with respect to a given probability density function $\rho:\mathbb{R}^m\rightarrow\mathbb{R}_+$,
\begin{equation}
\int f(\vx)^2\,\rho(\vx)\,d\vx \;<\; \infty,
\end{equation}
where we assume the domain of $f$ is the support of $\rho$. Given $\mU\in\mathbb{R}^{m\times n}$ with $n<m$ and $g:\mathbb{R}^n\rightarrow\mathbb{R}$, we measure the error in the ridge approximation with the $L^2(\rho)$ norm,
\begin{equation}
\label{eq:norm}
\left\| f(\vx) - g(\mU^T\vx) \right\|_{L^2(\rho)}
\;=\;
\left(
\int (f(\vx) - g(\mU^T\vx))^2\,\rho(\vx)\,d\vx
\right)^{\frac{1}{2}}.
\end{equation}
We restrict attention to matrices $\mU$ with orthonormal columns, $\mU^T\mU=\mI$, where $\mI$ is the $n\times n$ identity matrix. For a more general matrix with full column rank, we can transform to the orthonormal column case with a thin QR factorization, where the R factor represents an invertible change of variables in $g$'s domain. 

Given $\mU$ with orthonormal columns, let $\mV$ be an orthogonal basis for the complement of $\span(\mU)$ in $\mathbb{R}^m$, where $\span(\mU)$ denotes the span of $\mU$'s columns. The density function $\rho(\vx)$ induces joint, marginal, and conditional densities on the subspace coordinates of $\span(\mU)$ and $\span(\mV)$ as follows. For $\vy\in\mathbb{R}^n$ and $\vz\in\mathbb{R}^{m-n}$, define the following:
\begin{equation}
\begin{array}{rcll}
\pi(\vy,\vz) &=& \rho(\UyVz) & \text{(joint density)}\\
\pi(\vy) &=& \int \pi(\vy,\vz)\,d\vz & \text{(marginal density)}\\
\pzgy &=& \pi(\vy,\vz)/\pi(\vy) & \text{(conditional density)}
\end{array}
\end{equation}
The conditional density enables construction of a particularly useful ridge approximation. Define the conditional average of $f$ given subspace coordinates $\vy$, denoted $\mu$, as 
\begin{equation}
\label{eq:mu}
\mu 
\;=\; 
\mu(\vy, \mU) 
\;=\; 
\int f(\UyVz)\,\pzgy\,d\vz.
\end{equation}
Consider the ridge function $\mu(\mU^T\vx, \mU)$. By construction,
\begin{equation}
\label{eq:samemean}
\int \left(
\mu(\vy) - f(\UyVz)
\right)\,\pzgy\,d\vz \;=\; 0,
\end{equation}
for all $\vy$ such that $\pi(\vy)>0$. 
%Observe that the average of $f$ is the average of $\mu$,
%\begin{equation}
%\begin{aligned}
%\int \mu(\mU^T\vx, \mU)\,\rho(\vx)\,d\vx 
%&=
%\int \left( \int \mu(\vy,\mU) \pzgy\,d\vz \right)\,\py\,d\vy\\
%&=
%\int \mu(\vy,\mU) \,\py\,d\vy\\
%&=
%\int \left(\int f(\UyVz)\,\pzgy\,d\vz \right) \,\py\,d\vy\\
%&=
%\int f(\vx)\,\rho(\vx)\,d\vx,
%\end{aligned}
%\end{equation}
%which implies
%\begin{equation}
%\label{eq:samemean}
%\int f(\vx) - \mu(\mU^T\vx,\mU)\,\rho(\vx)\,d\vx \;=\; 0.
%\end{equation}
As a consequence of Pinkus's Theorem 8.3~\cite{pinkus2015}, for fixed $\mU$, \eqref{eq:samemean} implies that $\mu(\mU^T\vx,\mU)$ is the unique best ridge approximation in the $L^2(\rho)$ norm; see the discussion immediately following the theorem's statement. 

The particular choice of basis $\mU$ does not affect the ridge approximation $\mu$. In other words, we can replace $\mU$ by $\mU\mQ$, where $\mQ$ is an $n\times n$ orthogonal rotation matrix, and $\mu$ does not change. To see this, first examine the conditional density,
\begin{equation}
\begin{aligned}
\pi(\vz|\vy=\mQ^T\mU^T\vx)
&= 
\frac{\rho(\mU\mQ\mQ^T\mU^T\vx + \mV\vz)}{\int \rho(\mU\mQ\mQ^T\mU^T\vx + \mV\vz)\,d\vz}\\
&= 
\frac{\rho(\mU\mU^T\vx + \mV\vz)}{\int \rho(\mU\mU^T\vx + \mV\vz)\,d\vz}\\
&=
\pi(\vz|\vy=\mU^T\vx).
\end{aligned}
\end{equation}
Next examine the definition of $\mu$,
\begin{equation}
\begin{aligned}
\mu(\mQ^T\mU^T\vx;\,\mU\mQ)
&=
\int f(\mU\mQ\mQ^T\mU^T\vx + \mV\vz)\,\pi(\vz|\vy=\mQ^T\mU^T\vx)\,d\vz\\
&= 
\int f(\mU\mU^T\vx + \mV\vz)\,\pi(\vz|\vy=\mU^T\vx)\,d\vz\\
&=
\mu(\mU^T\vx;\,\mU).
\end{aligned}
\end{equation}
This implies that $\mu$ only depends on the subspace $\span(\mU)$ as opposed to the particular basis. For the rest of this section, the notation $\mU$ denotes an equivalence class of matrices whose columns span the same subspace. Similarly, we use $\mV$ to represent an equivalence class of matrices whose columns span the orthogonal complement of $\span(\mU)$ in $\mathbb{R}^m$.

To characterize the optimal $\mU$, we derive a differentiable cost function from \eqref{eq:norm}. Define $R=R(\mU)$ as
\begin{equation}
\label{eq:R}
R(\mU) \;=\; \frac{1}{2}\,\left\|
f(\vx)-\mu(\mU^T\vx,\mU)
\right\|_{L^2(\rho)}^2.
\end{equation}
Note that, similar to $\mu$, $R$ only depends on $\span(\mU)$ as opposed to the choice of basis. Therefore, the appropriate manifold for optimization is the Grassmann manifold---i.e., the space of $n$-dimensional subspaces of $\mathbb{R}^m$, denoted $\mathbb{G}(n,m)$. Let $\mU_\ast$ be a solution to the following optimization problem:
\begin{equation}
\label{eq:min}
\begin{array}{ll}
\underset{\mU}{\operatorname{minimize}
} & R(\mU),\\
\text{subject to} & \mU\in\mathbb{G}(n,m).
\end{array}
\end{equation}
We call $\mU_\ast$ an \emph{optimal subspace}. In general, the objective function is not a convex function of $\mU$, so its minimizer may not be unique. In practice, we use numerical methods to estimate $\mU_\ast$. 

It is convenient to reformulate the optimization \eqref{eq:min} in terms of the complement subspace $\span(\mV)$. The conditional average $\mu$ in \eqref{eq:mu} is the average of $f(\vx)$ over the affine subspace $S(\vx)$ defined as
\begin{equation}
S(\vx) \;=\; \{\,
\vx'\in\mathbb{R}^m \,\mid\, \vx'=\mU\mU^T\vx + \mV\vz,\,\vz\in\mathbb{R}^{m-n}
\,\}.
\end{equation}
This space depends only on the shift $\mU\mU^T\vx$ and $\span(\mV)$---not the choice of basis for $\span(\mV)$. We can write the shift as
\begin{equation}
\mU\mU^T\vx \;=\; (\mI -\mV\mV^T) \,\vx,
\end{equation}
where $\mI$ is the $m\times m$ identity matrix. Again, this shift does not depend on the choice of basis---only the subspace $\span(\mV)$. Therefore, we can write $\mu$ from \eqref{eq:mu} as
\begin{equation}
\label{eq:mu2}
\mu
\;=\;
\mu(\vx, \mV)
\;=\;
\int f\left( (\mI-\mV\mV^T)\vx + \mV\vz\right)\,\pzgy\,d\vz.
\end{equation}
Similarly, we rewrite $R$ from \eqref{eq:R} as
\begin{equation}
R(\mV) \;=\; \frac{1}{2}\,\left\| 
f(\vx) - \mu(\vx, \mV) 
\right\|^2_{L^2(\rho)}.
\end{equation}
Let $\mV_\ast$ be an $(m-n)$-dimensional subspace that satisfies
\begin{equation}
\label{eq:minV}
\begin{array}{ll}
\underset{\mV}{\operatorname{minimize}
} & R(\mV),\\
\text{subject to} & \mV\in\mathbb{G}(m-n,m).
\end{array}
\end{equation}
An optimal $\mU_\ast$ that solves \eqref{eq:min} is the orthogonal complement of a particular $\mV_\ast$. 

Reformulating $R$ as a function of $\mV$ is convenient for studying its gradient. Edelman et al.~\cite[Section 2.5.3]{Edelman1998} derive a formula for the gradient of $R$ on the Grassmann manifold in terms of the partial derivatives on the ambient Euclidean space $\mathbb{R}^{m\times(m-n)}$. Denote the gradient on the Grassmann by $\bnabla$. Then
\begin{equation}
\label{eq:RV}
\bnabla R(\mV) 
\;=\; 
\ppV R(\mV) - \mV\mV^T \ppV R(\mV)
\;=\;
\mU\mU^T \ppV R(\mV),
\end{equation}
where $\ppV R$ is the $m\times (m-n)$ matrix of partial derivatives
\begin{equation}
\left(\ppV R\right)_{ij} \;=\; \frac{\partial R}{\partial v_{ij}},\qquad i=1,\dots,m,\quad j=1,\dots,m-n,
\end{equation}
where $v_{ij}$ is the $(i,j)$ element of $\mV$. This formula can be implemented and passed to a gradient-based nonlinear optimization routine, e.g., steepest descent or a quasi-Newton method. 

%Numerical optimization methods are sensitive to the initial guess for the optimizer; for example, Newton's method converges quadratically to a local minimum \emph{when the starting point is sufficiently close to the local minimum}. 

%% file: sec2-subspaces.tex
\section{A near-stationary subspace}
\label{sec:stationary}

\noindent The objective function $R(\mV)$ in \eqref{eq:minV} is, in general, not a convex function of $\mV$, so a gradient-based optimization algorithm is only guaranteed to converge to a stationary point\footnote{Recent work suggests that the probability of terminating at a stationary point that is not a local minimizer is zero~\cite{lee2016}.}. Additionally, the cost of reaching a stationary point may depend heavily on the initial guess. We call a subspace \emph{near-stationary} if we can bound the norm of the objective's gradient at that subspace.

\begin{definition}
\label{def:stationary}
A subspace $\mV_\ast\in\mathbb{G}(m-n,m)$ is near-stationary if there is a constant $A=A(f,\rho)$ such that
\begin{equation}
\left\| \bnabla R(\mV_\ast) \right\|_F
\;\leq\; A,
\end{equation}
where $\|\cdot\|_F$ is the Frobenius norm, and $\bnabla R$ is the gradient on the Grassmann manifold of $R$ from \eqref{eq:minV}.
\end{definition}

In what follows, we show that the active subspace derived from $f$'s derivatives is near-stationary. The active subspace is the eigenspace of a particular symmetric, positive semidefinite matrix, and the bound $A$ from Definition \ref{def:stationary} is related to the matrix's eigenvalues. In statistical regression, Samarov~\cite{Samarov1993} studied related matrices built from derivatives of the regression's link function, which he termed \emph{average derivative functionals}; Samarov's $\mT_1$ is similar to the matrix we study. The regression case contrasts ours since we assume $f$ and its derivatives are given, whereas the link function in regression depends on parameters to be estimated from data.

To ensure that all necessary quantities exist, we make the following assumption on $f(\vx)$.

\begin{assumption}
\label{assump1}
Given the probability density $\rho:\mathbb{R}^m\rightarrow\mathbb{R}_+$, assume (i) $f\in L^2(\rho)$ is continuous and differentiable for all $\vx$ in the support of $\rho$ and (ii) the partial derivatives of $f$ are continuous and square-integrable with respect to $\rho$,
\begin{equation}
\int \left(\frac{\partial f}{\partial x_i}(\vx)\right)^2\,\rho(\vx)\,d\vx \;<\; \infty.
\end{equation}
\end{assumption}
For $f$ that satisfies Assumption \ref{assump1}, define the $m\times m$ symmetric positive semidefinite matrix $\mC=\mC(f,\rho)$ as
\begin{equation}
\label{eq:C}
\mC \;=\; \int \nabla f(\vx)\,\nabla f(\vx)^T\,\rho(\vx)\,d\vx
\;=\; \mW\Lambda\mW^T,
\end{equation}
where $\Lambda$ is the diagonal matrix of nonnegative eigenvalues $\lambda_1,\dots,\lambda_m$ in decreasing order, and $\mW$ is the orthogonal matrix of corresponding eigenvectors. Assume also that $\lambda_n>\lambda_{n+1}$ for some $n<m$, and consider the partition
\begin{equation}
\label{eq:eigsep}
\Lambda = \bmat{\Lambda_1 & \\ & \Lambda_2},\qquad
\mW = \bmat{\mW_1 & \mW_2},
\end{equation}
where $\Lambda_1$ contains the first $n$ eigenvalues, and $\mW_1$ contains the first $n$ eigenvectors. The active subspace~\cite{Constantine2014,asmbook} is the span of the columns of $\mW_1$. The eigenvalues reveal whether $f$ is a ridge function, as seen in the following theorem.

\begin{theorem}
\label{thm:necsuff}
For $f$ that satisfies Assumption \ref{assump1}, assume that $\lambda_n>\lambda_{n+1}$ for some $n<m$. A vector $\vw$ is in the null space of $\mC$ from \eqref{eq:C} if and only if $f(\vx)$ is constant along $\vw$ for all $\vx$ in the support of $\rho$.
\end{theorem}
The proof of Theorem \ref{thm:necsuff} is in \ref{sec:app1}. We next consider a ridge approximation constructed with the subspace $\span(\mW_1)$; the next theorem, which is Theorem 3.1 from~\cite{Constantine2014}, bounds the $L^2(\rho)$ approximation error. The key to the following theorem is a weighted Poincar\'{e} inequality for the weight function $\rho$. Edmunds and Opic~\cite{Edmunds1993} provide details of the Poincar\'{e} constant with weight function $\rho$, and Bebendorf~\cite{Bebendorf2003} proves a Poincar\'{e} inequality for convex domains. Chen~\cite{chen1982inequality} shows that the Poincar\'{e} constant is 1 when $\rho$ is a standard normal density, and he proves a version of the Poincar\'{e} inequality for a correlated normal density. The following theorem assumes that $\rho$ admits a Poincar\'{e} inequality. 

\begin{theorem}
\label{thm:approx}
For $f$ that satisfies Assumption \ref{assump1} and $\rho$ that admits a Poincar\'{e} inequality, define $\mu$ as in \eqref{eq:mu}. Then 
\begin{equation}
\left\|f(\vx) -\mu(\mW_1^T\vx,\mW_1)\right\|_{L^2(\rho)} \;\leq\;
C\,(\lambda_{n+1} + \cdots + \lambda_m)^{\frac{1}{2}},
\end{equation}
where $C=C(\rho)$ is the Poincar\'{e} constant associated with the probability density function $\rho$. 
\end{theorem}
The proof of Theorem \ref{thm:approx} is in Section 3 of~\cite{Constantine2014} and Section 4.2 of~\cite{asmbook}. Note that if we write $\mu$ as $\mu(\vx,\mV)$, as in \eqref{eq:mu2}, then the Theorem \ref{thm:approx} holds for $\mu(\vx,\mW_2)$. The next theorem shows that the active subspace is near-stationary when $\rho(\vx)$ is a standard Gaussian density and $f(\vx)$ is Lipschitz continuous.

\begin{theorem}
\label{thm:gradbased}
Let $\rho$ be a standard Gaussian density on $\mathbb{R}^m$, and assume that $f$ satisfies Assumption \ref{assump1} with $\rho$ a Gaussian density. Additionally, assume that 
\begin{itemize}
\item[(i)] $f$ is Lipschitz continuous with constant $L$,
\item[(ii)] $\lambda_n>\lambda_{n+1}$ for some $n<m$.
\end{itemize}
Then for $R$ as in \eqref{eq:minV} and $\mW_2$ from \eqref{eq:eigsep},
\begin{equation}
\label{eq:gradbnd}
\left\|\bnabla R(\mW_2)\right\|_F
\;\leq\; 
L\,\left(
2m^{\frac{1}{2}} + (m-n)^{\frac{1}{2}}
\right)\,
(\lambda_{n+1} + \cdots + \lambda_m)^{\frac{1}{2}},
\end{equation}
where $\bnabla$ denotes the gradient on $\mathbb{G}(m-n,m)$, and $\|\cdot\|_F$ is the Frobenius norm.
\end{theorem}
The proof of Theorem \ref{thm:gradbased} is in \ref{sec:app2}. The bound's dependence on the eigenvalues implies that if $f$ is a ridge function of $n$ variables, then $\span(\mW_1)$ is a stationary point for the minimization \eqref{eq:min}. We expect that the Gaussian assumption on $\rho$ can be relaxed at the cost of a more complicated bound in \eqref{eq:gradbnd} involving the gradient of $\rho$. Such an extension is beyond the scope of this manuscript.

%% file: sec3-examples.tex
\section{Computational examples}
\label{sec:examples}

\noindent Theorems \ref{thm:necsuff} and \ref{thm:gradbased} suggest a computational heuristic for fitting a ridge approximation. Assuming the gradient $\nabla f(\vx)$ can be evaluated as a subroutine (e.g., via algorithmic differentiation~\cite{Griewank00}), consider the steps in Algorithm \ref{alg:asridge}.

\begin{algorithm}
\begin{enumerate}
\itemsep0em
\item Estimate the matrix $\mC$ from \eqref{eq:C} with a numerical integration rule, and compute its eigendecomposition.
\item Choose $n$ such that $\lambda_n>\lambda_{n+1}$ and $\lambda_{n+1},\dots,\lambda_m$ are relatively small.
\item Use the first $n$ estimated eigenvectors as an initial guess for numerical optimization of \eqref{eq:min}.
\end{enumerate}
\caption{Exploiting the active subspace for ridge approximation}
\label{alg:asridge}
\end{algorithm}

Constantine and Gleich~\cite{constantine2015computing} analyze a Monte Carlo method for estimating $\mC$ and its eigendecomposition as in step 1 of Algorithm \ref{alg:asridge}. If the estimated eigenvalues do not decay appropriately to choose $n$ in step 2, then the given $f(\vx)$ may not be a good candidate for ridge approximation. It is easy to construct functions that are not amenable to ridge approximation, e.g., $f(\vx)=\|\vx\|^2$ or any radially symmetric function; such structure would manifest as little-to-no decay in the eigenvalues. In Section \ref{sec:goodex}, we offer a computational heuristic for step 3 of Algorithm \ref{alg:asridge} based on alternating minimization. 

\subsection{A simple example where the heuristic fails}
\label{sec:badex}

\noindent The heuristic in Algorithm \ref{alg:asridge} relies on $\mC$'s eigenvalues to measure the suitability of the associated eigenvectors for an initial guess when fitting a ridge approximation. We show a bivariate example where there is a large gap between the first and second eigenvalues, but the second eigenvector---though a stationary point---is far from the global minimizer of the objective function \eqref{eq:R}. In the bivariate case, we can parameterize the rotation in the two-dimensional domain by one angle $\alpha\in[0,\pi]$.

Let $\rho(x_1,x_2)$ be a standard bivariate Gaussian density, and consider the bivariate function
\begin{equation}
\label{eq:bv}
f(x_1,x_2) \;=\; 5x_1 + \sin(10\pi x_2).
\end{equation}
This function has a Lipschitz constant $L$ that is bounded by 32. The matrix $\mC$ from \eqref{eq:C} is (to 4 significant digits)
\begin{equation}
\mC \;=\; \bmat{25.00 & 0 \\ 0 & 526.4}
\;=\;
\underbrace{\bmat{0 & 1 \\ 1 & 0}}_{\mW}\underbrace{\bmat{526.4 & 0 \\ 0 & 25.00}}_{\Lambda}\bmat{0 & 1 \\ 1 & 0}^T.
\end{equation}
We estimate $\mC$ with a tensor product Gauss-Hermite quadrature rule with 101 points per dimension (10201 total points), which was sufficient for four digits of accuracy. 

\begin{figure}[!h]
\centering
\includegraphics[width=0.5\textwidth]{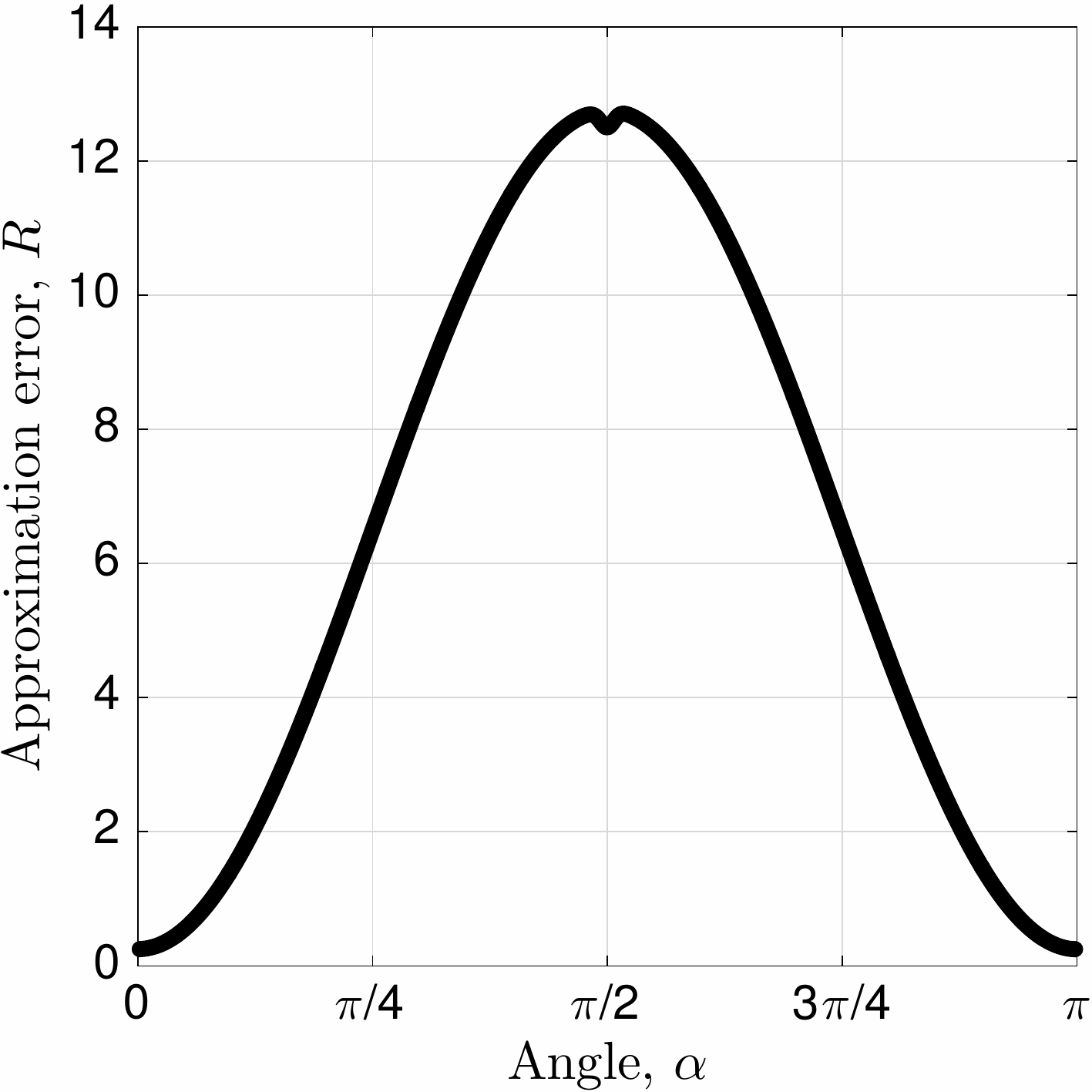}%
\caption{The $L^2(\rho)$ error \eqref{eq:R} as a function of subspace angle $\alpha$ for the ridge approximation of $f(x_1,x_2)=5x_1+\sin(10\pi x_2)$. The active subspace is $\span([0,1]^T)$---corresponding to $\alpha=\pi/2$---which is a poor initial guess for a gradient-based optimizer.}
\label{fig:alpha}
\end{figure}

According to the heuristic in Algorithm \ref{alg:asridge}, the eigenvalues $\Lambda$ suggest that the vector $[0,1]^T$, which corresponds to $\alpha=\pi/2$, would be a good starting point for a numerical optimization. Figure \ref{fig:alpha} plots the error $R$ as a function of the subspace angle $\alpha$ for 500 values of $\alpha\in[0,\pi]$. Each $R$ is computed with Gauss-Hermite quadrature rule with 301 points in each dimension (90601 total points), which is sufficient for four digits of accuracy. The figure shows that $[0,1]^T$ (i.e., $\alpha=\pi/2$) is actually a local minimizer of $R$ with $R([0,1]^T)=12.5$. A gradient-based optimization routine starting at $[0,1]^T$ is unlikely to escape the local minimum. In contrast, $R([1,0]^T)=0.25$, where $\span([1,0]^T)$ (corresponding to $\alpha=0$) is the orthogonal complement of the active subspace. In other words, the eigenvector associated with the smaller eigenvalue is both a stationary point and a minimizer. 

This example suggests a type of function for which the heuristic is not well suited, namely, functions that oscillate rapidly along one direction and vary slowly but consistently along another. The derivative-based metrics choose the direction of oscillation as the important direction even when a ridge approximation is more accurate, in the mean-squared sense, along another direction. 

\subsection{A ridge approximation response surface for drag in a transonic airfoil}
\label{sec:goodex}

%\begin{algorithm}[H]
%\label{alg:alt}
%\KwData{A data set $\{\vx_i,f(\vx_i)\}$, an initial subspace $\mU_0$, and a polynomial degree $N$.}
%\KwResult{Estimate the optimal subspace $\mU_\ast$ and optimal polynomial parameters $\vp_\ast$.}
%\For{$i=1$ to maxiter}{
%Compute $\vy_i = \mU_0^T\vx_i$.\;
%Use least-squares to fit parameters $\vp_\ast$ of polynomial approximation $p_N(\vy,\vp)$ with data $\{\vy_i,f(\vx_i)\}$.\;
%Compute $\mU_\ast$ as the minimizer of 
%}
%\caption{Alternating minimization.}
%\end{algorithm}

\noindent The following numerical example demonstrates the suitability of the heuristic in Algorithm \ref{alg:asridge} applied to a model of drag coefficient as a function of airfoil shape in a standard transonic test problem; more details on the model follow. We emphasize that the goal of the study is to assess the quality of the gradient-based active subspaces as an initial guess for a numerical optimization method for the ridge approximation problem. It is not our goal to assess the quality of the ridge approximation; however, we do present metrics on the approximation quality for completeness. 

First, we propose a computational heuristic for estimating the minimizer of \eqref{eq:min} based on a polynomial model and alternating minimization~\cite{Ruhe1980}; this corresponds to step 3 of Algorithm \ref{alg:asridge}. The polynomial model plays the role of $\mu$ in the approximation error \eqref{eq:R}. Let $p_N(\vy,\theta)$ be a polynomial of degree $N$ in $n$ variables (i.e., $\vy\in\mathbb{R}^n$), where $\theta$ is the vector of the polynomial's parameters (i.e., coefficients). The dimension of $\theta$ depends on the number of terms in the polynomial model, which depends on $N$ and $n$. In our experiments, we use a multivariate polynomial model of total degree $N$, so the number of terms is ${N+n \choose n}$~\cite{Dunkl2014}. This model is more general than the projection pursuit regression model \eqref{eq:ppr}, since it includes products of powers of the $n$ linear combinations. 

%Assume $M$ input/output pairs $(\vx_i,f_i)$ are given, where $f_i=f(\vx_i)$. Choose an initial matrix $\mU_0\in\mathbb{R}^{m\times n}$ and a degree $N$ for the polynomial approximation $p_N$. Consider the following iteration:

\begin{algorithm}
Given $M$ input/output pairs $(\vx_i,f(\vx_i))$, $\mU_0\in\mathbb{R}^{m\times n}$ with orthogonal columns, polynomial degree $N$, and number of iterations $P$.
\vskip 1em
\noindent For $i$ from 1 to $P$, do
\begin{enumerate}
\item Compute $\vy_i=\mU_0^T\vx_i$ for $i=1,\dots,M$.
\item Compute $\theta_\ast$ as the solution to the least-squares problem,
\begin{equation}
\label{eq:polyapprox}
\begin{array}{ll}
\underset{\theta}{\operatorname{minimize}} & \sum_{i=1}^M \left(f_i - p_N(\vy_i,\theta)\right)^2.
\end{array}
\end{equation}
\item Compute $\mU_\ast$ as the solution to the Grassmann manifold-constrained least-squares problem,
\begin{equation}
\label{eq:polyresid}
\begin{array}{ll}
\underset{\mU}{\operatorname{minimize}} & \sum_{i=1}^M \left(f_i - p_N(\mU^T\vx_i,\theta_\ast)\right)^2,\\
\text{subject to} & \mU\in\mathbb{G}(n,m).
\end{array}
\end{equation}
\item Set $\mU_0=\mU_\ast$.
\end{enumerate}
\caption{Polynomial-based alternating minimization scheme for \eqref{eq:min}}
\label{alg:alternate}
\end{algorithm}

Algorithm \ref{alg:alternate} warrants several comments. We use an alternating scheme over the two sets of variables, $\theta$ and $\mU$, because of its simplicity. Alternating schemes are known to stall and/or converge very slowly relative to gradient-based approaches using all variables~\cite[Section 9.3]{optbook}. For this reason, we do not offer specific stopping criteria in Algorithm \ref{alg:alternate}. Instead, we opt for a user-defined number $P$ of iterations, which requires more intervention from the user; in the experiment below, we use $P=20$, which was sufficient to demonstrate the efficiency of the active subspace as a starting point $\mU_0$. Also, the gradient of the objective function in \eqref{eq:polyresid} is much easier to implement than $R$ from \eqref{eq:min} or \eqref{eq:minV}, since $p_N(\vy,\theta_\ast)$ is independent of $\mU$---unlike $\mu=\mu(\vy,\mU)$ from \eqref{eq:mu}. However, analyzing the ridge approximation with $p_N$ is much more difficult. The Python codes that implement Algorithm \ref{alg:alternate} can be found at \url{bitbucket.org/paulcon/near-stationary-subspace}. We use the package Pymanopt~\cite{Pymanopt2016} to estimate the Grassmann manifold-constrained least-squares problem in \eqref{eq:polyresid} with the gradient-based steepest descent method.

Relative to standard polynomial-based response surfaces, the ridge approximation can---for the same number $M$ of function evaluations $(\vx_i,f_i)$---fit a higher degree polynomial along the directions that $f(\vx)$ varies. In other words, with $M$ fixed in \eqref{eq:polyapprox}, the degree $N$ can be much larger in $n$ variables than in $m>n$ variables. However, if several iterations of the alternating heuristic are needed to achieve the stopping criteria, then fitting the ridge approximation may itself be costly due to the relatively expensive Grassmann-constrained minimization step. Therefore, a good initial subspace can be very advantageous, as seen in the following example. 

We apply the alternating minimization heuristic to build a ridge approximation response surface for an aerospace engineering model of a transonic airfoil's drag coefficient as a function of its shape; details on this model are in~\cite[Section 5.3]{asmbook}. The baseline airfoil is the NACA0012 (a standard transonic test case for computational fluid dynamics), and perturbations to the baseline shape are parameterized by $m=18$ Hicks-Henne bump functions. Each of the 18 parameters is constrained to the interval $[-0.01,0.01]$ to ensure valid airfoil geometries, and we choose $\rho(\vx)$ to be a uniform density on the hypercube $[-0.01,0.01]^{18}$. Note that this density does not satisfy the Gaussian assumption of Theorem \ref{thm:gradbased}, but this does not impede us from numerically testing the heuristic in Algorithm \ref{alg:asridge}. Given the airfoil geometry, the drag coefficient is computed with the Stanford University Unstructured (SU2) compressible Euler solver~\cite{su2}. This software also solves the continuous adjoint equation for the Euler equations, which enables the computation of the gradient of the drag coefficient with respect to the 18 shape parameters. To summarize, $f(\vx)$ is the airfoil's drag coefficient as a function of the shape parameters, and a computer model returns $f$ and $\nabla f$ given $\vx$. Preliminary tests with the model (e.g., standard grid convergence tests and parameter sweeps) indicate that, using our chosen mesh and solver tolerances in SU2, we can expect at least four digits of accuracy in each evaluation of drag and its partial derivatives over the parameter space, which was sufficient for our experiments. One evaluation of both drag and its gradient takes approximately 2 minutes on one core; in total, we ran 20000 simulations for the following experiment. Given the data set of 20000 simulations, all computations for the following experiments were run on two nodes of the Colorado School of Mines Mio cluster. Each experiment used 4 cores, which accelerated the numpy operations. 

\begin{figure}[!h]
\centering
\subfloat[Eigenvalue error estimates]{
\label{fig:eigerr}
\includegraphics[width=0.48\textwidth]{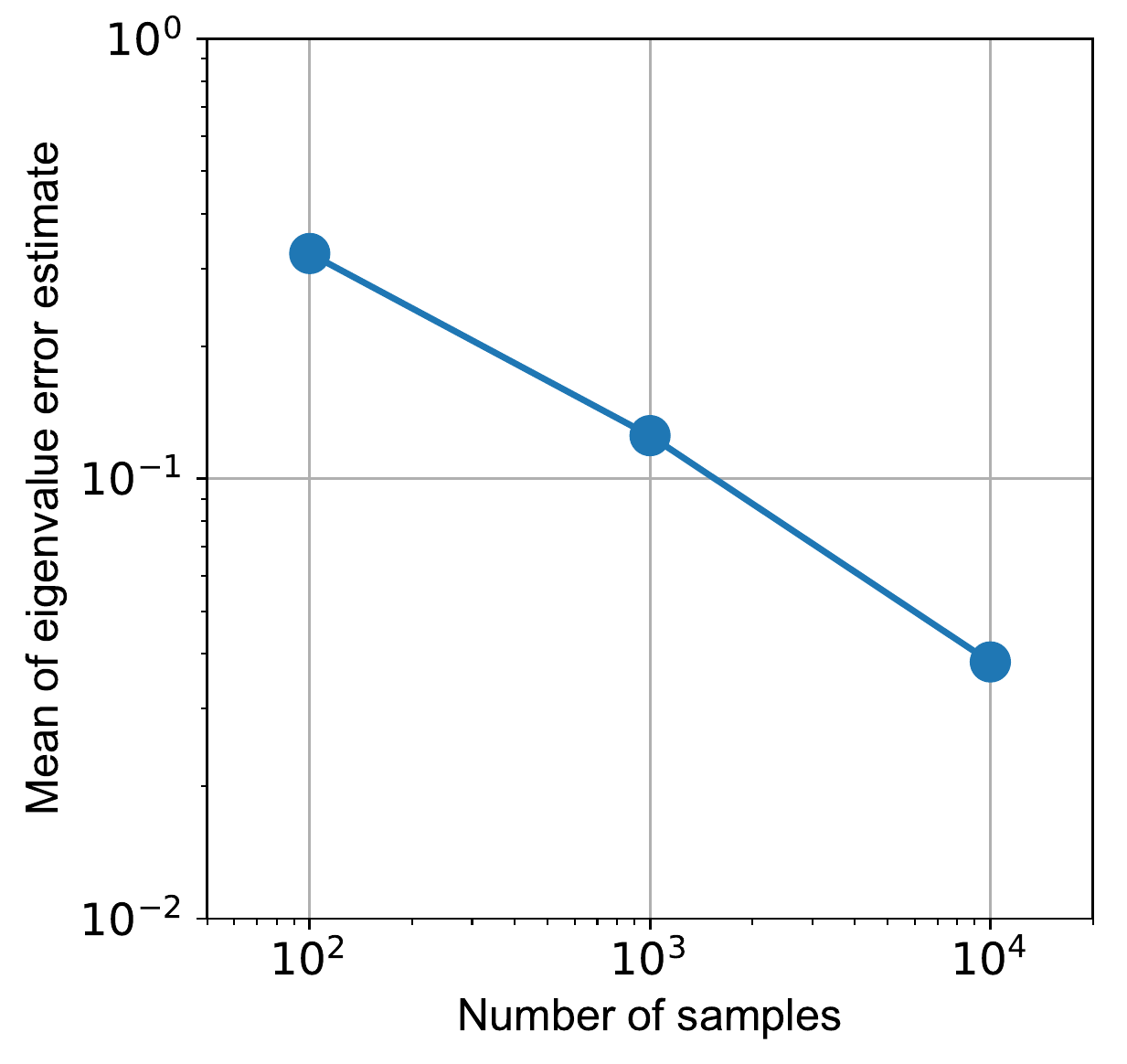}%
}
\hfil
\subfloat[Subspace error estimates]{
\label{fig:suberr}
\includegraphics[width=0.48\textwidth]{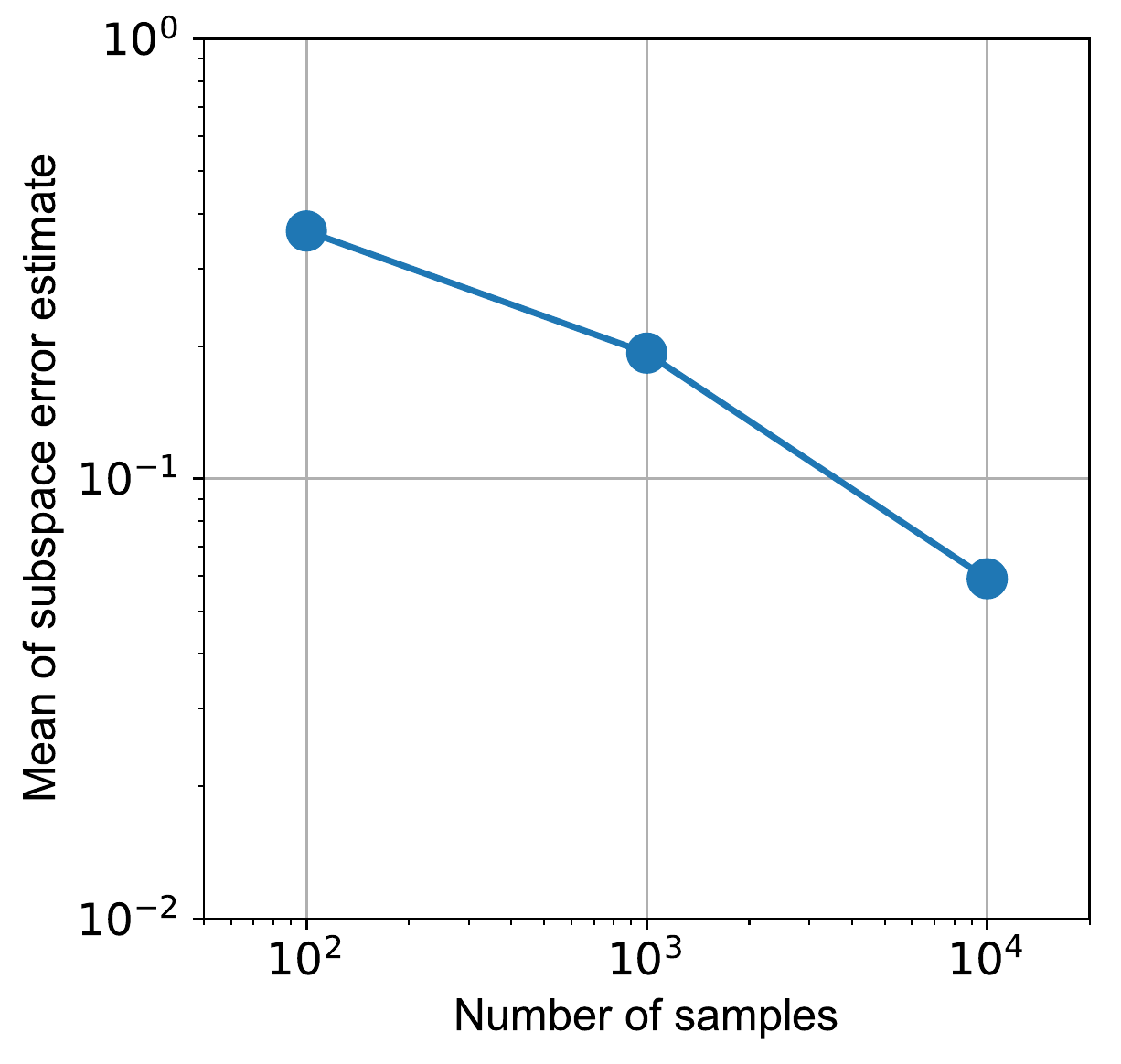}%
}
\caption{Eigenvalue and subspace errors \eqref{eq:errors} with respect to reference values computed from a 20000-sample Latin hypercube design.}
\label{fig:aserrs}
\end{figure}

\begin{figure}[!h]
\centering
\subfloat[Eigenvalue estimates]{
\label{fig:eigest}
\includegraphics[width=0.48\textwidth]{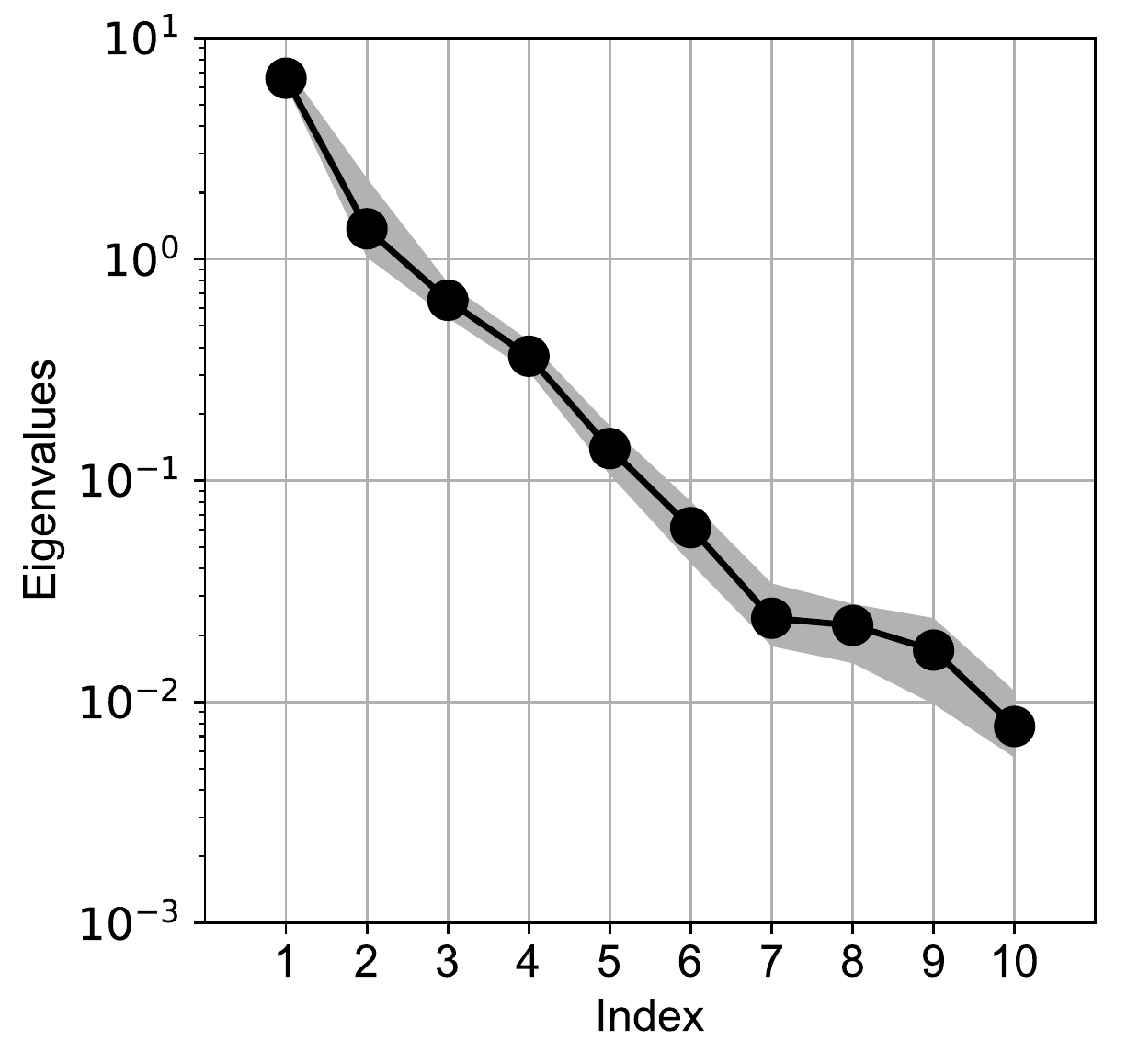}%
}
\hfil
\subfloat[Subspace distance estimates]{
\label{fig:suberrest}
\includegraphics[width=0.48\textwidth]{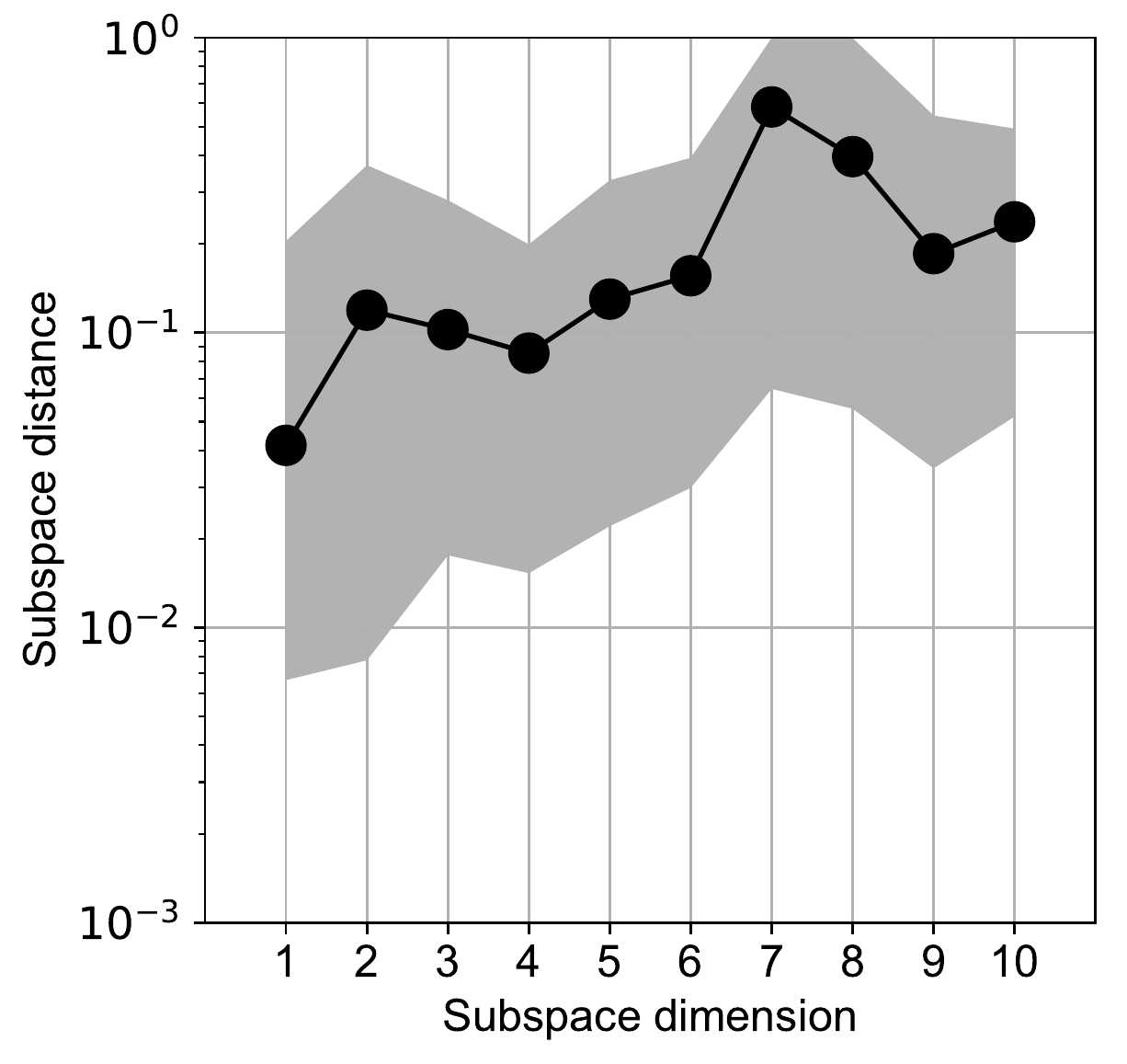}%
}
\caption{Bootstrap-based variability estimates for eigenvalues and subspaces using the first $N=1000$ samples from the 20000-sample Latin hypercube design.}
\label{fig:eigsnsub}
\end{figure}

We estimate $\mC$ from \eqref{eq:C} with Latin hypercube sampling~\cite{McKay1979}. Constantine and Gleich~\cite{constantine2015computing} analyzed a simple Monte Carlo method for estimating active subspaces. We first generate 20000 Latin hypercube samples $\{\vx_i\}$, and we compute drag $f(\vx_i)$ and its gradient vector $\nabla f(\vx_i)$ (computed via the adjoint equations in SU2) for each sample. We use all 20000 samples to compute reference values $\mCref = \mWref\Lambdaref\mWref^T$. We then select the first 100, 1000, and 10000 samples from the Latin hypercube design and estimate $\mC$ with the associated subsets of $\{\nabla f(\vx_i)\}$; note that the subsets of $\{\vx_i\}$ do not satisfy the Latin property of the original Latin hypercube design~\cite{Qian2009}. To verify the expected $\mathcal{O}(N^{-1/2})$ of the Monte Carlo scheme, we compute the following two error metrics:
\begin{equation}
\label{eq:errors}
\begin{aligned}
\text{err}_{\Lambda,N} &= \frac{1}{m}
\sum_{k=1}^m \frac{| \lambda_{\text{ref},k} - \lambda_{N,k} |}{| \lambda_{\text{ref},k} |}, \\
\text{err}_{\mW,N} &= \frac{1}{m-1}
\sum_{k=1}^{m-1} 
\left\| 
\mW_{\text{ref},k}\mW_{\text{ref},k}^T - \mW_{N,k}\mW_{N,k}^T 
\right\|_2,
\end{aligned}
\end{equation}
where a subscript $k$ on $\lambda$ indicates the $k$th eigenvalue (ordered in descending order), a subscript $k$ on $\mW$ indicates the first $k$ eigenvectors, a subscript $\text{ref}$ indicates the reference value, and the subscript $N$ indicates using the first $N$ samples, where $N\in\{100,1000,10000\}$. The error $\text{err}_{\Lambda,N}$ is the average relative error across eigenvalue estimates. Figure \ref{fig:eigerr} shows this error as a function of $N$. The convergence rate of $\mathcal{O}(N^{-1/2})$ appears as expected, which verifies the Monte Carlo implementation. The error $\text{err}_{\mW,N}$ is the average subspace error over subspaces of dimension 1 through $m-1$. Each term in the sum is the matrix 2-norm difference between projector matrices onto their respective subspaces, which is also the sine of the principal angle between subspaces; this is a common measure of subspace distance~\cite[Chapter 2]{gvl13}. Figure \ref{fig:suberr} shows this error as a function $N$, and it shows the expected convergence rate, which verifies the Monte Carlo implementation. 

For a fixed number of samples, Constantine and Gleich~\cite{constantine2015computing} propose a bootstrap-based heuristic for estimating the variability in the Monte Carlo estimates of the eigenvalues and subspaces, which is similar to the bootstrap-based eigenvalue standard error estimates from Efron and Tibshirani~\cite[Chapter 7.1]{efron1994}. Since the eigenvalue and subspace estimates are nonlinear functions of the data, a central limit theorem-based standard error is not possible. For a data set with $N$ elements, the bootstrap method draws $N$ samples with replacement from the data; computing the desired quantities from this data set (e.g., the eigenvalues of $\mC$'s estimate) yields a \emph{bootstrap replicate}. Repeating this procedure many times produces a set of bootstrap replicates---eigenvalues and subspaces, in our case---for the same original data set of $N$ elements. 

Figure \ref{fig:eigest} shows (i) the first 10 eigenvalue estimates as black dots and (ii) the bootstrap ranges (min/max over the replicates) in the shaded region for the first $N=1000$ samples from the 20000-sample Latin hypercube design. The tight bootstrap ranges suggest that there is not much variability in the eigenvalue estimates under perturbations to the data. For each bootstrap replicate of the first $k$ eigenvectors, $\mW_{N,k}^\ast$, the $k$-dimensional subspace distance is computed as 
\begin{equation}
\label{eq:bssuberr}
\text{err}_{N,k}^\ast \;=\; 
\left\| 
\mW_{N,k}\mW_{N,k}^T - \left(\mW_{N,k}^\ast\right)\left(\mW_{N,k}^\ast\right)^T
\right\|_2.
\end{equation}
Figure \ref{fig:suberrest} shows the mean (black dots) and min/max range (shaded region) over these subspace distance replicates as a function of subspace dimension. We emphasize that calculations represented by Figures \ref{fig:eigest} and \ref{fig:suberrest} do not use reference values; they only include a subset of $N=1000$ samples from the 20000-sample Latin hypercube design. 

The bootstrap-based subspace variability estimates behave consistently with the error analysis from Gleich and Constantine~\cite{constantine2015computing}. In particular, a small subspace error for subspace dimension $n$ corresponds to a large gap between eigenvalues $\lambda_n$ and $\lambda_{n+1}$; this is consistent with well known perturbation results for invariant subspaces~\cite{Stewart1973}. Constantine et al.~\cite{Constantine2014} showed that errors in the active subspace have substantial impact on errors in the ridge approximation $\mu(\mW_1^T\vx,\mW_1)$ from Theorem \ref{thm:gradbased}. However, for an optimization-based ridge approximation from Algorithm \ref{alg:alternate}, since the numerically estimated active subspace is used only as an initial guess, we expect errors in the subspace to matter little; this is supported by the following experiment.

\begin{figure}[!h]
\centering
\subfloat[$n=1$, $N=2$]{
\includegraphics[width=0.24\textwidth]{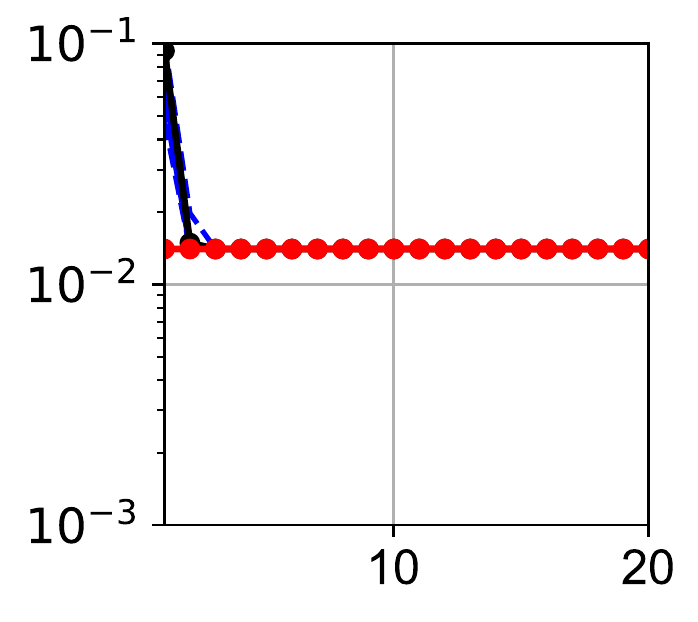}%
}
\subfloat[$n=1$, $N=3$]{
\includegraphics[width=0.24\textwidth]{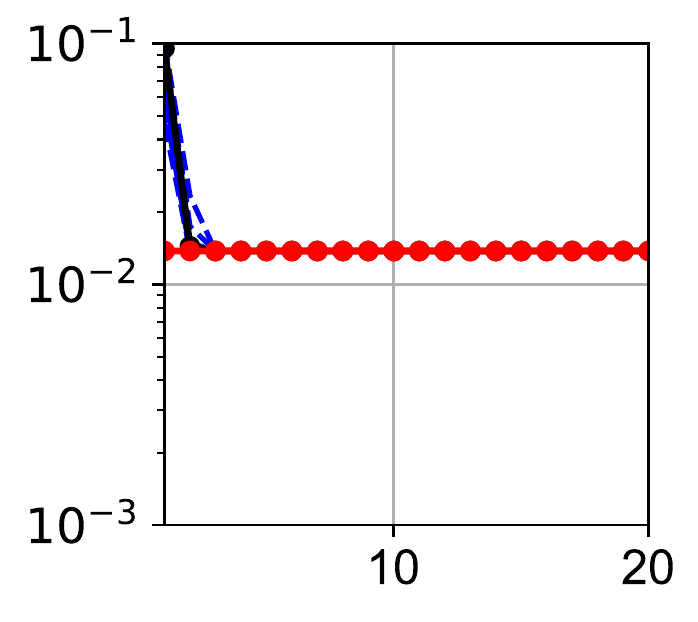}%
}
\subfloat[$n=1$, $N=4$]{
\includegraphics[width=0.24\textwidth]{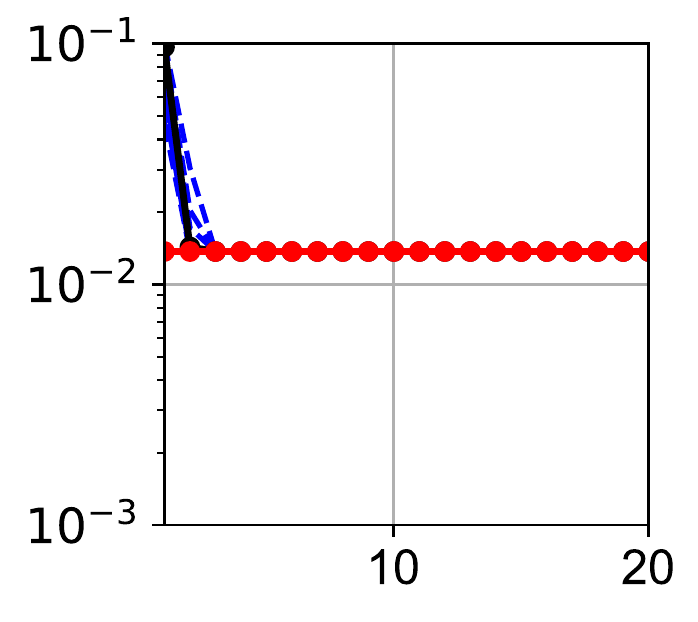}%
}
\subfloat[$n=1$, $N=5$]{
\includegraphics[width=0.24\textwidth]{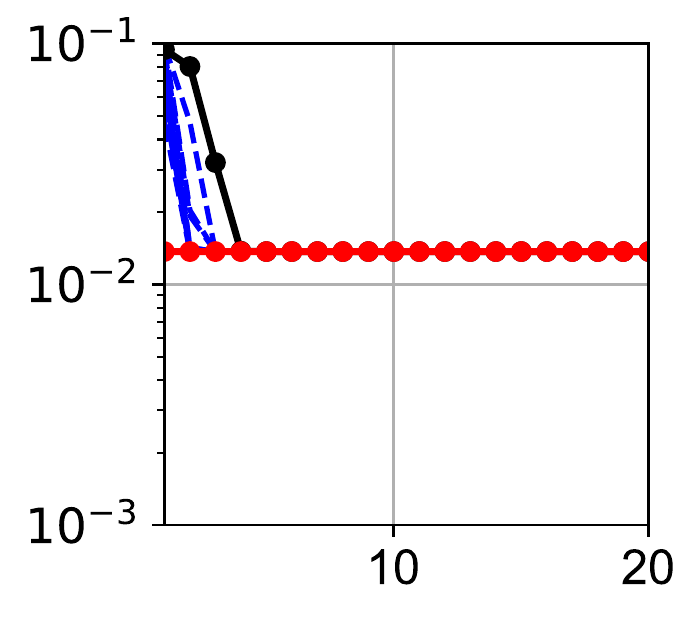}%
}
\\
\subfloat[$n=2$, $N=2$]{
\includegraphics[width=0.24\textwidth]{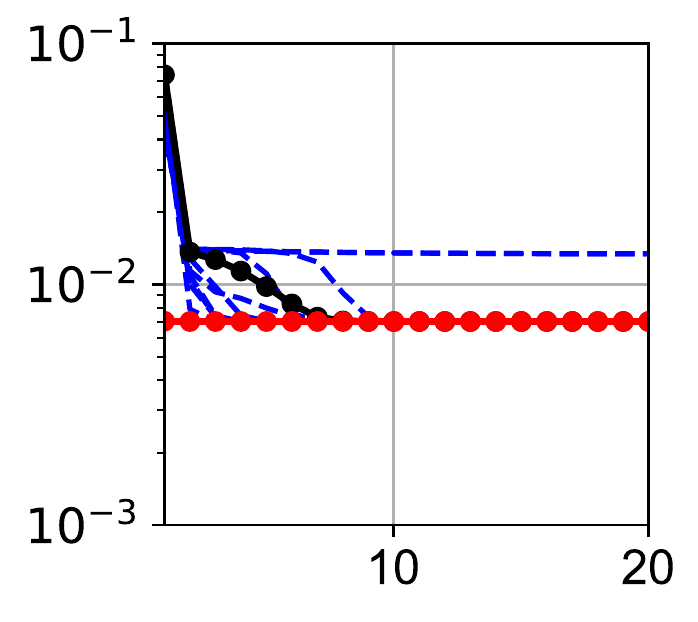}%
}
\subfloat[$n=2$, $N=3$]{
\includegraphics[width=0.24\textwidth]{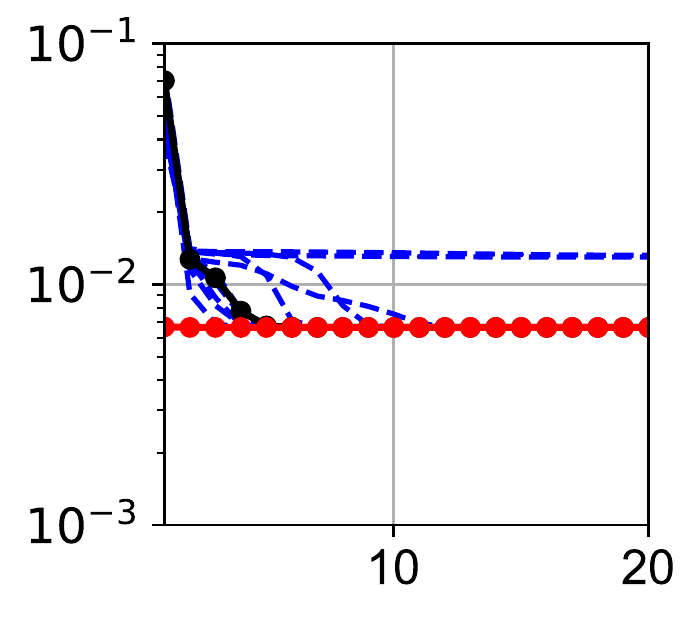}%
}
\subfloat[$n=2$, $N=4$]{
\includegraphics[width=0.24\textwidth]{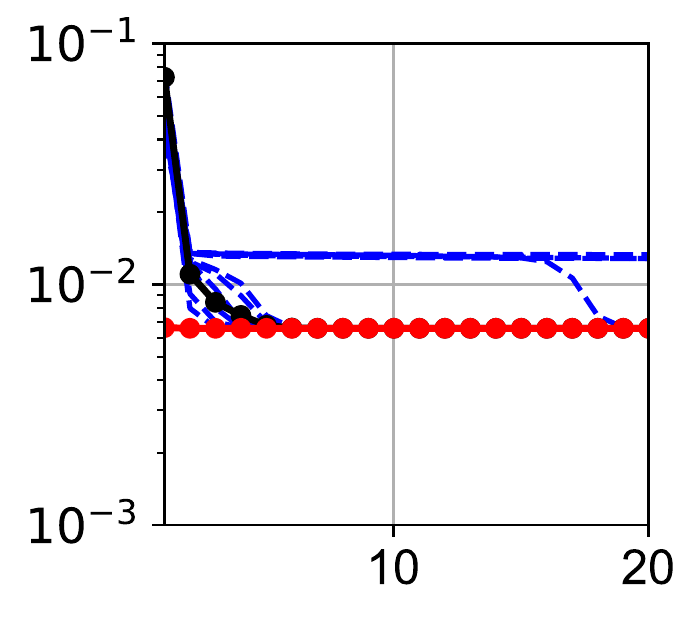}%
}
\subfloat[$n=2$, $N=5$]{
\includegraphics[width=0.24\textwidth]{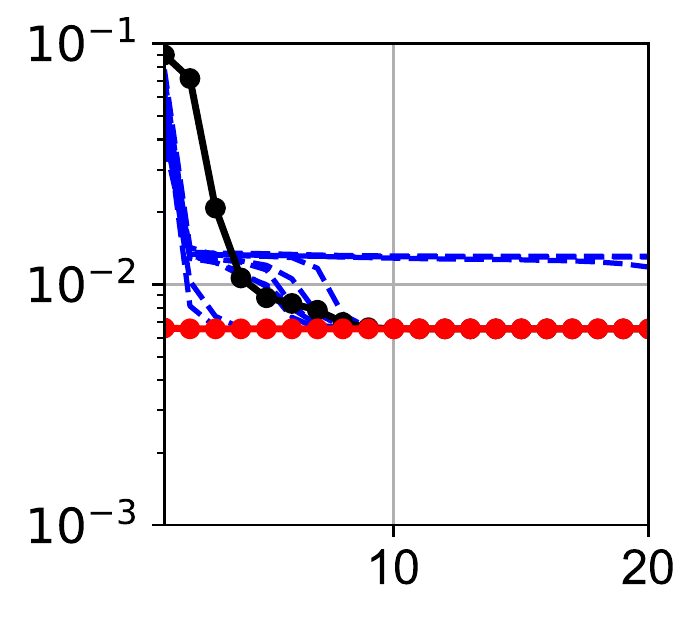}%
}
\\
\subfloat[$n=3$, $N=2$]{
\includegraphics[width=0.24\textwidth]{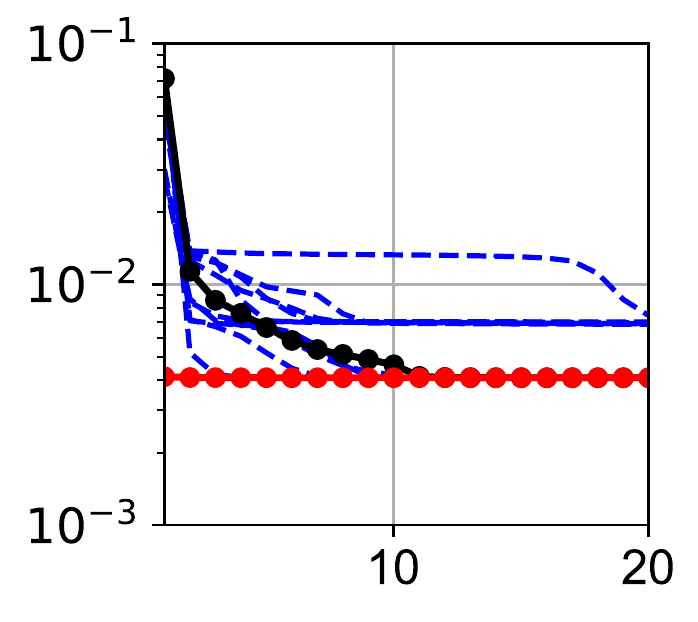}%
}
\subfloat[$n=3$, $N=3$]{
\includegraphics[width=0.24\textwidth]{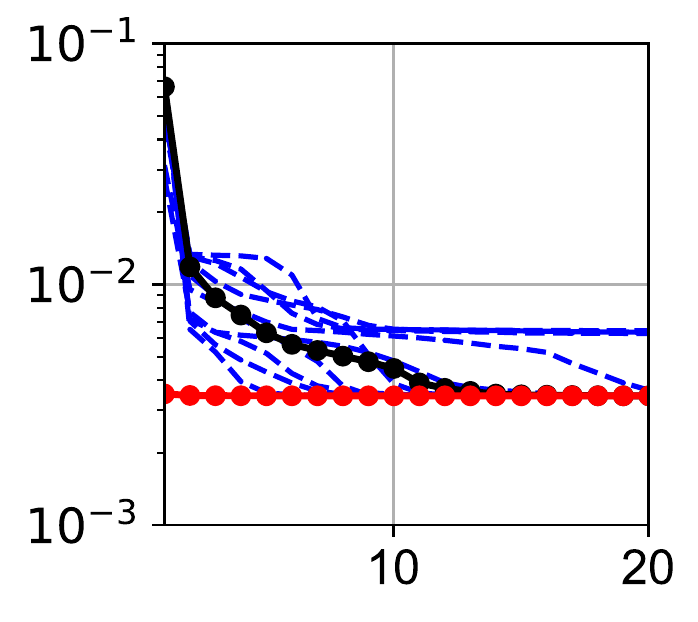}%
}
\subfloat[$n=3$, $N=4$]{
\includegraphics[width=0.24\textwidth]{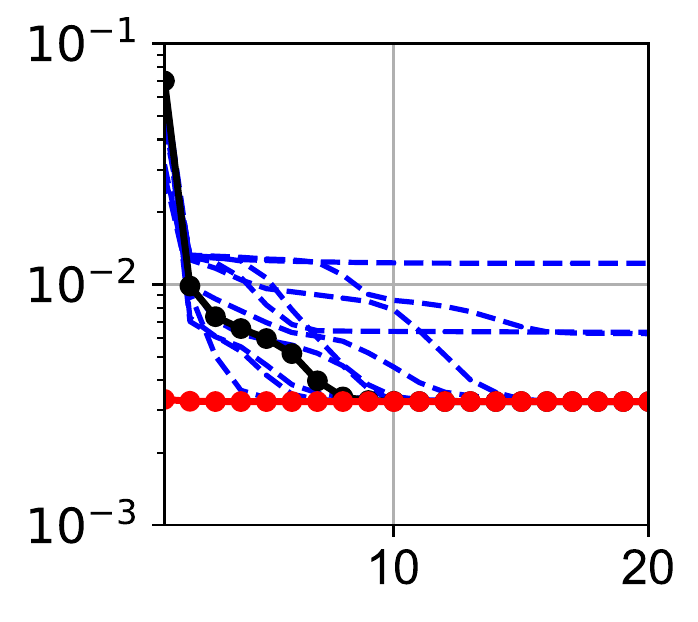}%
}
\subfloat[$n=3$, $N=5$]{
\includegraphics[width=0.24\textwidth]{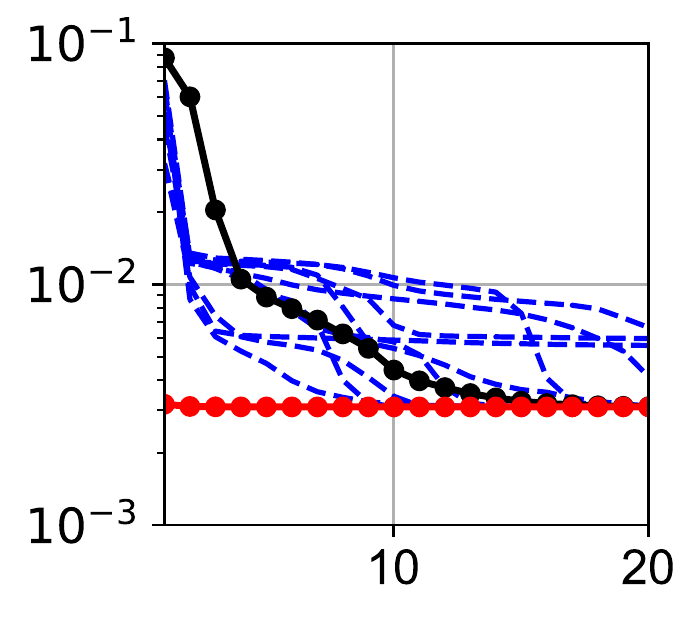}%
}
\\
\subfloat[$n=4$, $N=2$]{
\includegraphics[width=0.24\textwidth]{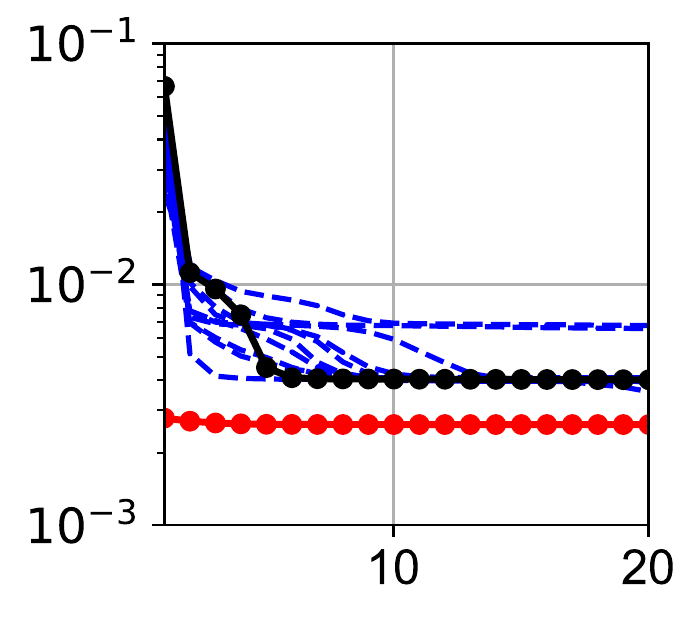}%
}
\subfloat[$n=4$, $N=3$]{
\includegraphics[width=0.24\textwidth]{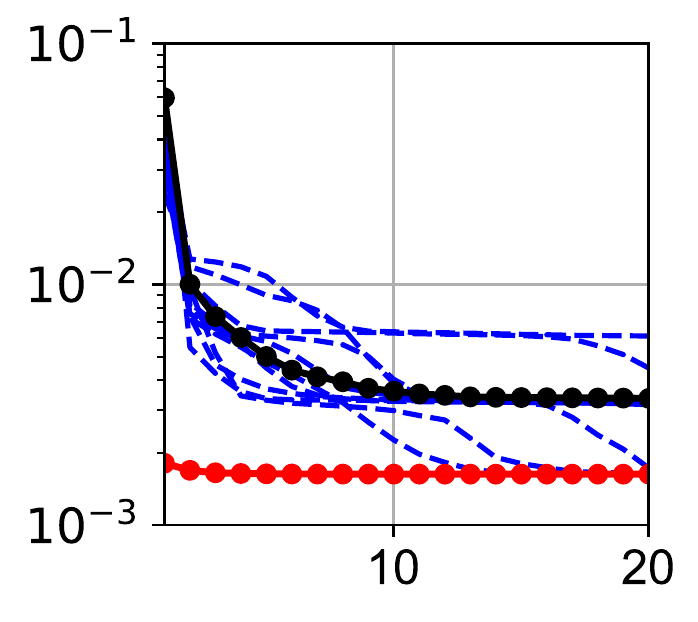}%
}
\subfloat[$n=4$, $N=4$]{
\includegraphics[width=0.24\textwidth]{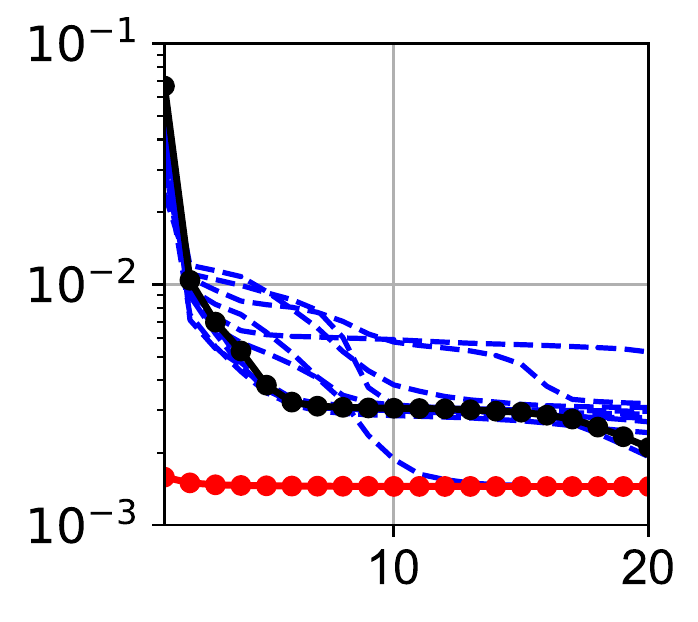}%
}
\subfloat[$n=4$, $N=5$]{
\includegraphics[width=0.24\textwidth]{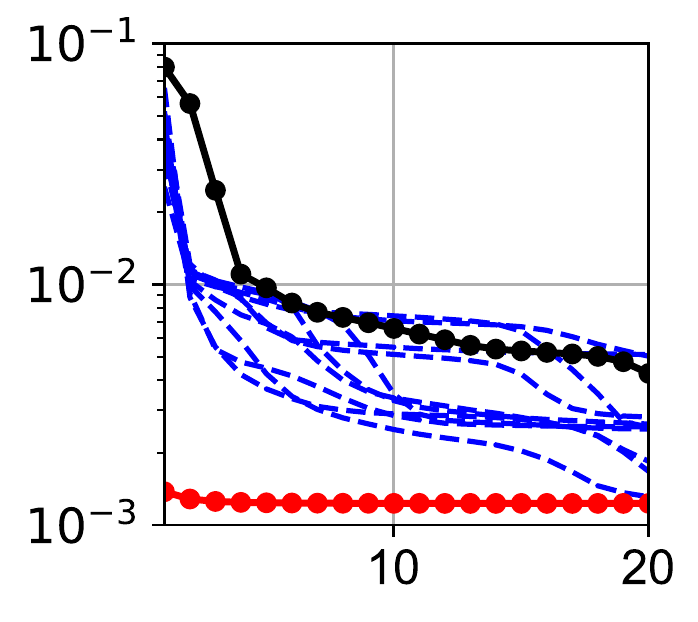}%
}
\caption{Each subfigure shows the residual \eqref{eq:polyresid}---using the first 1000 samples from the 20000-sample Latin hypercube design---as a function of the iteration in the alternating minimization heuristic Algorithm \ref{alg:alternate}. The black connected dots use the first $n$ columns of the $m\times m$ identity matrix as $\mU_0$. The blue dashed lines show results using 10 random starting points for $\mU_0$. The red connected dots use the first $n$ eigenvectors of $\mCref$. The subfigures vary the number $n$ of linear combinations from 1 to 4 (top to bottom) and the degree $N$ of the polynomial approximation from 2 to 5 (left to right). In all cases, the $n$ eigenvectors of the numerically estimated $\mC$ from \eqref{eq:C} provide a superior starting point.}
\label{fig:iters}
\end{figure}

Figure \ref{fig:iters} shows the results of an experiment comparing different starting points $\mU_0$ for the alternating heuristic in Algorithm \ref{alg:alternate}: (i) a random $m\times n$ matrix with orthogonal columns, (ii) the first $n$ columns of the $m\times m$ identity matrix, and (iii) the first $n$ estimated eigenvectors from the 20000-sample reference estimate $\mCref$ of $\mC$ from \eqref{eq:C}. The data in each case is the first $M=1000$ pairs $(\vx_i, f(\vx_i))$---i.e., shape parameters and associated drag---from the 20000-sample Latin hypercube design of experiments. We ran similar experiments for the first 100, 500, 5000, and 10000 pairs---all results were qualitatively similar. Each subfigure in Figure \ref{fig:iters} shows the residual \eqref{eq:polyresid} as a function of the iteration count; the residual value is on the vertical axis, and the number of iterations is on the horizontal axis. For the Grassmann manifold optimization in \eqref{eq:polyresid}, we set the maximum number of steepest descent steps to 10, which made each iteration considerably faster than converging to a specified tolerance on the norm of the residual gradient on the Grassmann. The black connected dots show the results using the identity matrix starting point; the blue dashed lines show results from 10 different random starting points; and the red connected dots show results using the first $n$ eigenvectors of $\mCref$. The subfigures vary the polynomial degree $N$ from 2 to 5 (left to right) and the number $n$ of linear combinations from 1 to 4 (top to bottom). For every case, the first $n$ eigenvectors of $\mCref$ provide a superior starting point for the alternating heuristic, and the advantage increases as $N$ and $n$ increase. 

\begin{figure}[!h]
\centering
\subfloat[$n=1$, $N=2$]{
\includegraphics[width=0.24\textwidth]{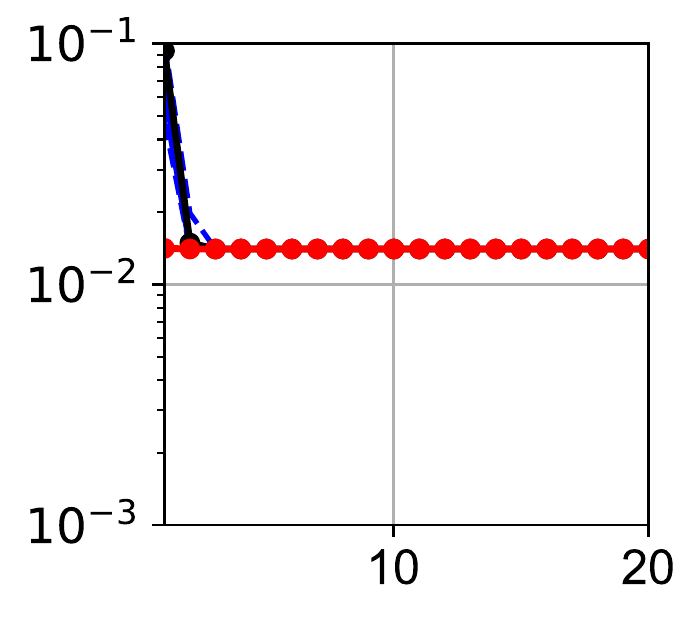}%
}
\subfloat[$n=1$, $N=3$]{
\includegraphics[width=0.24\textwidth]{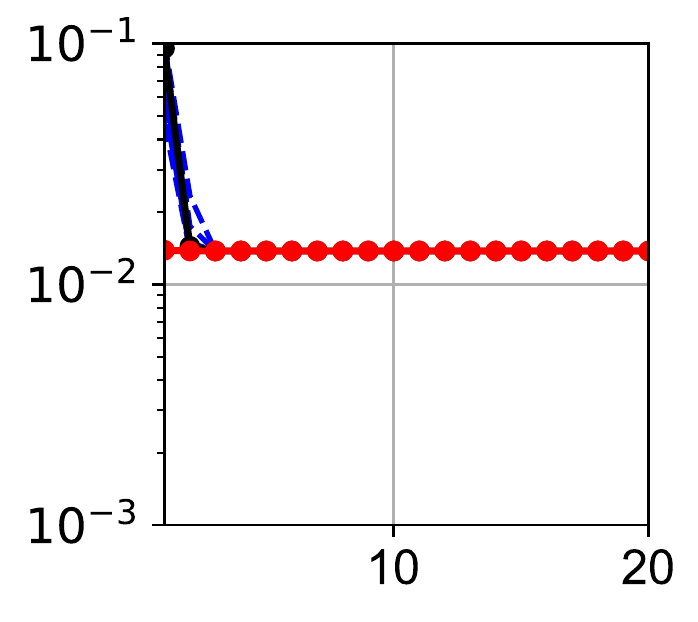}%
}
\subfloat[$n=1$, $N=4$]{
\includegraphics[width=0.24\textwidth]{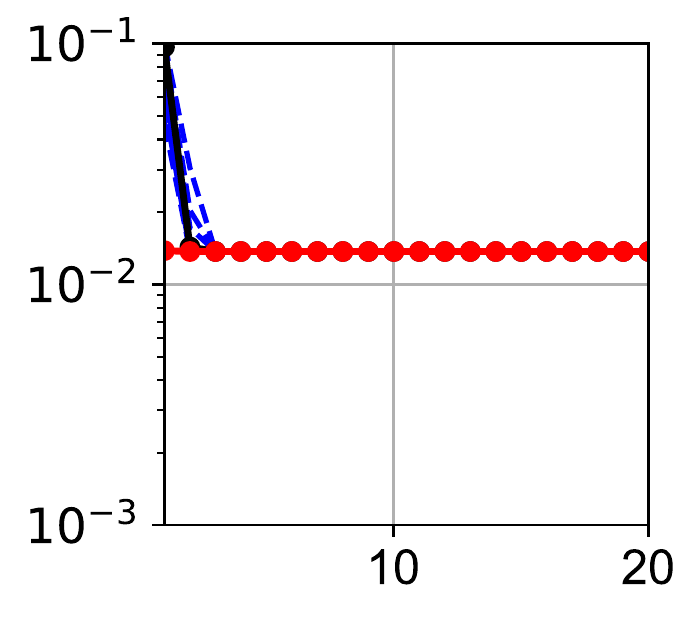}%
}
\subfloat[$n=1$, $N=5$]{
\includegraphics[width=0.24\textwidth]{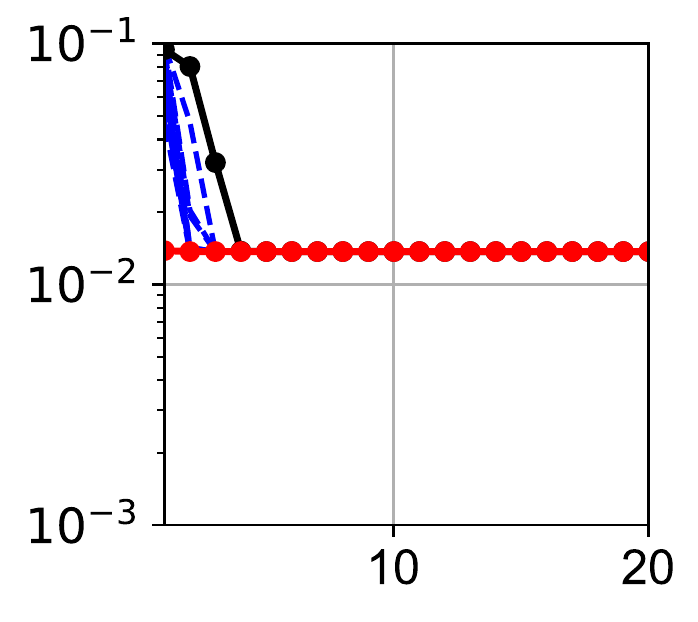}%
}
\\
\subfloat[$n=2$, $N=2$]{
\includegraphics[width=0.24\textwidth]{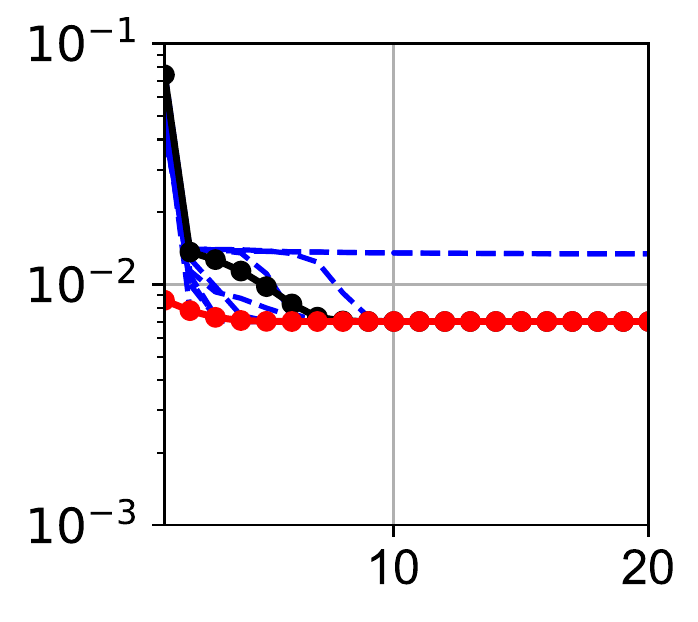}%
}
\subfloat[$n=2$, $N=3$]{
\includegraphics[width=0.24\textwidth]{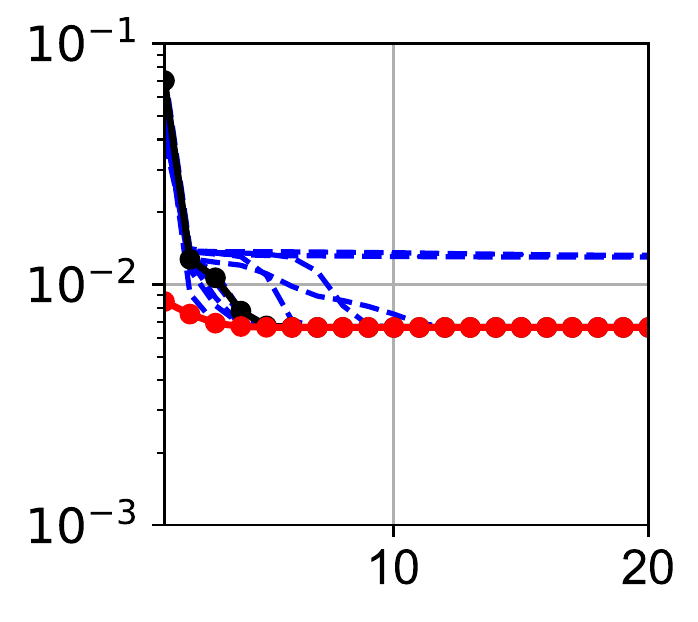}%
}
\subfloat[$n=2$, $N=4$]{
\includegraphics[width=0.24\textwidth]{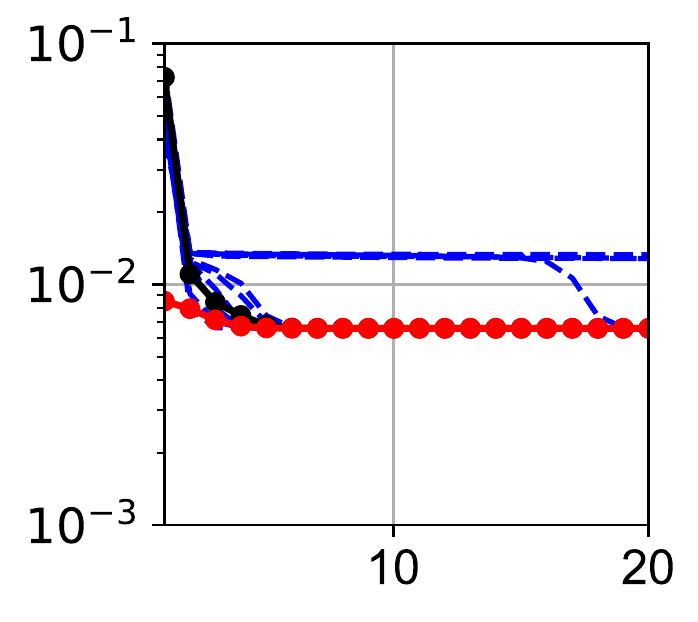}%
}
\subfloat[$n=2$, $N=5$]{
\includegraphics[width=0.24\textwidth]{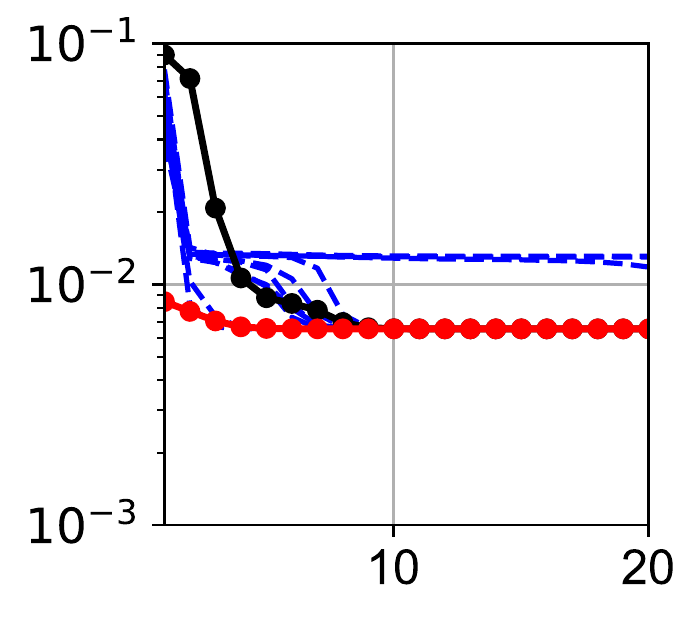}%
}
\\
\subfloat[$n=3$, $N=2$]{
\includegraphics[width=0.24\textwidth]{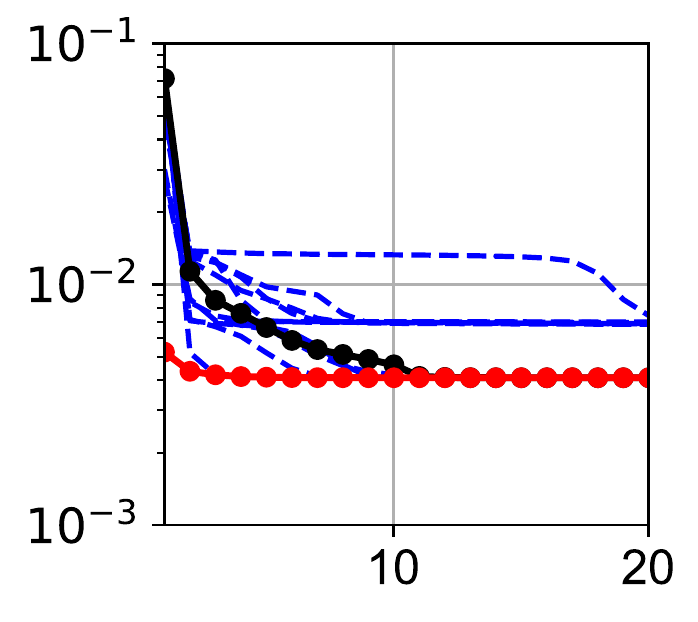}%
}
\subfloat[$n=3$, $N=3$]{
\includegraphics[width=0.24\textwidth]{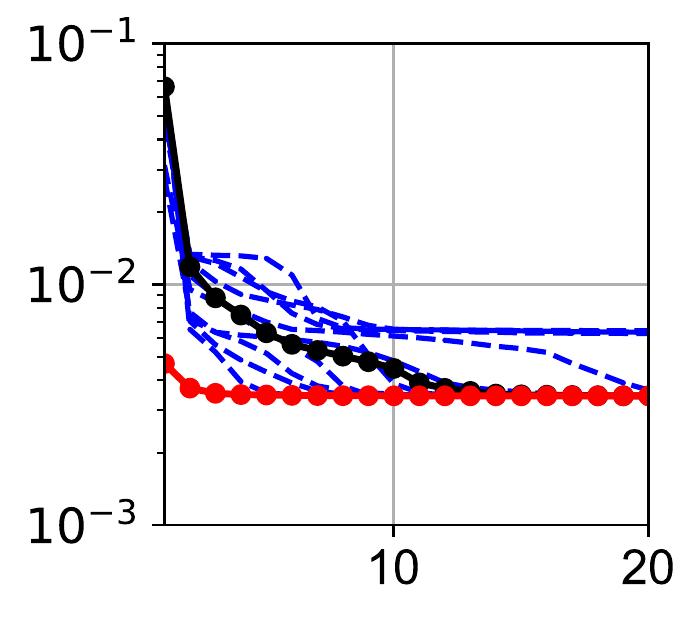}%
}
\subfloat[$n=3$, $N=4$]{
\includegraphics[width=0.24\textwidth]{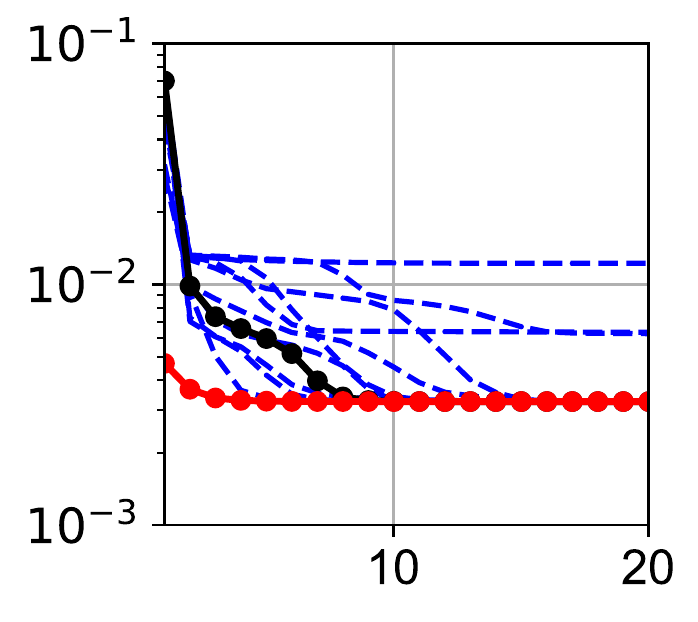}%
}
\subfloat[$n=3$, $N=5$]{
\includegraphics[width=0.24\textwidth]{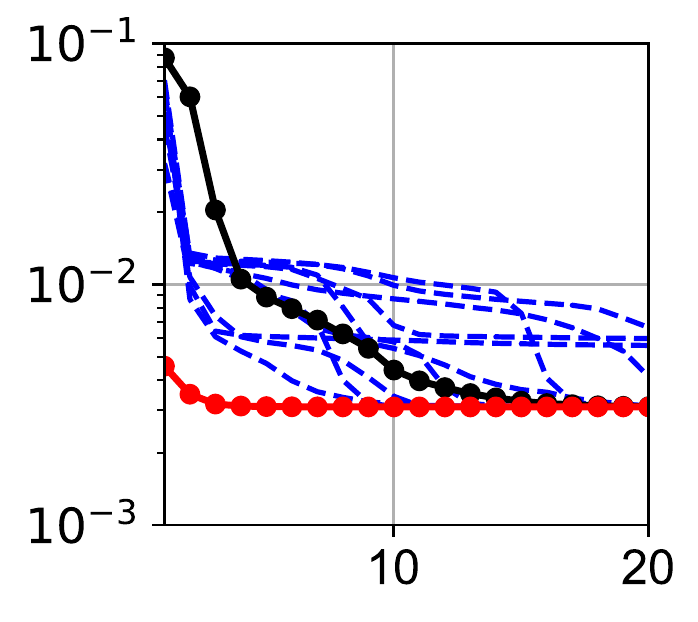}%
}
\\
\subfloat[$n=4$, $N=2$]{
\includegraphics[width=0.24\textwidth]{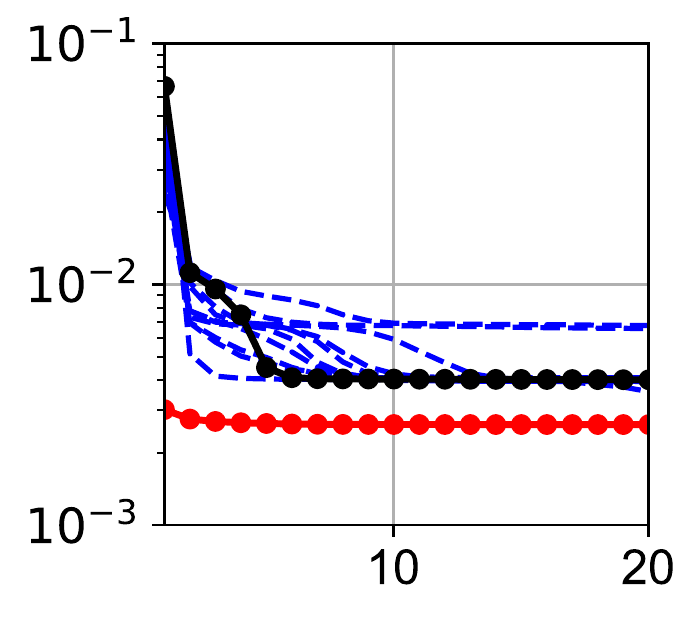}%
}
\subfloat[$n=4$, $N=3$]{
\includegraphics[width=0.24\textwidth]{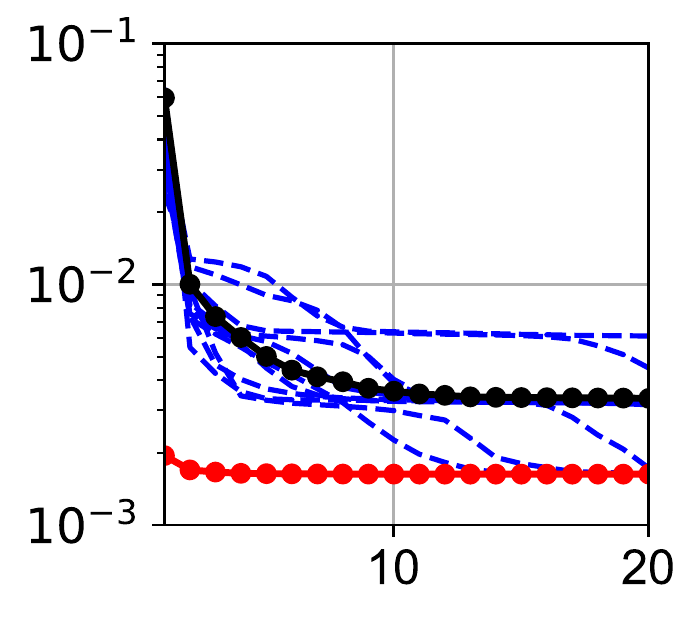}%
}
\subfloat[$n=4$, $N=4$]{
\includegraphics[width=0.24\textwidth]{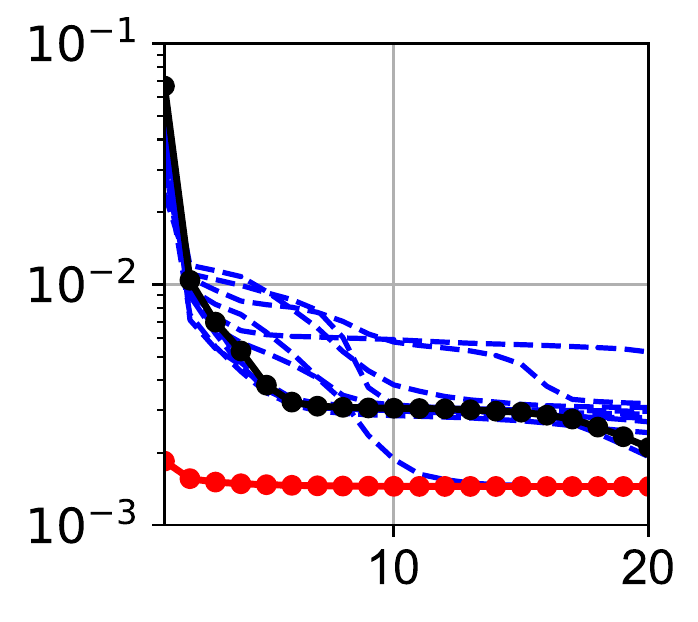}%
}
\subfloat[$n=4$, $N=5$]{
\includegraphics[width=0.24\textwidth]{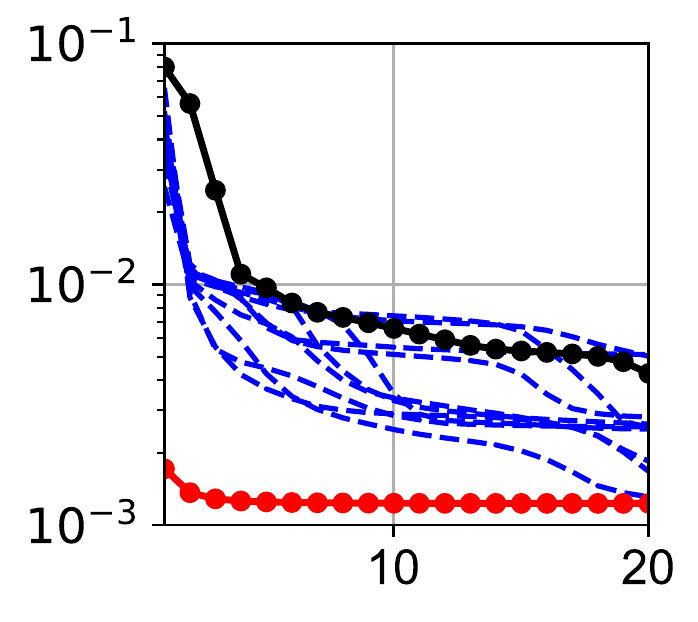}%
}
\caption{Identical to Figure \ref{fig:iters}, except the active subspace initial guesses are computed with only 10 gradient samples. Notice the slightly greater decrease in the first few iterations, compared to Figure \ref{fig:iters}, that results from the lower accuracy approximation.}
\label{fig:iters-sm}
\end{figure}

Computing the active subspace may take significant computational effort, whereas generating a random starting point or a few columns of the identity matrix takes relatively no computational effort. The experiment from Figure \ref{fig:iters} uses estimated eigenvectors from the very expensive 20000-sample reference value $\mCref$. In Figure \ref{fig:iters-sm}, we repeat the experiment using a 10-sample estimate of $\mC$'s eigenvectors; the results are remarkably similar. There is a marginally greater decrease in the residual from the initial cheap, 10-sample estimate of the active subspace relative to the 20000-sample reference value. This suggests that a cheap estimate of the active subspace may still be a better starting point than a random starting point for the ridge approximation heuristic in Algorithm \ref{alg:alternate}. For this particular data set of $M=1000$ runs, a marginal cost of 10 runs (with adjoint-based gradients) to compute a good initial subspace for ridge approximation is quite reasonable. The specific trade-offs for other models and data sets will depend on the relative cost of estimating the eigenvectors of $\mC$ versus fitting the ridge approximation.

\begin{figure}[!h]
\centering
\subfloat[$n=1$, $N=2$]{
\includegraphics[width=0.24\textwidth]{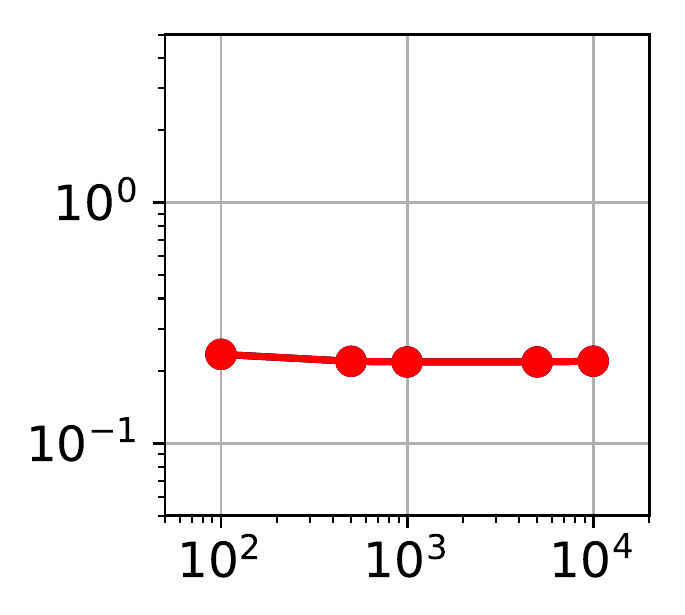}%
}
\subfloat[$n=1$, $N=3$]{
\includegraphics[width=0.24\textwidth]{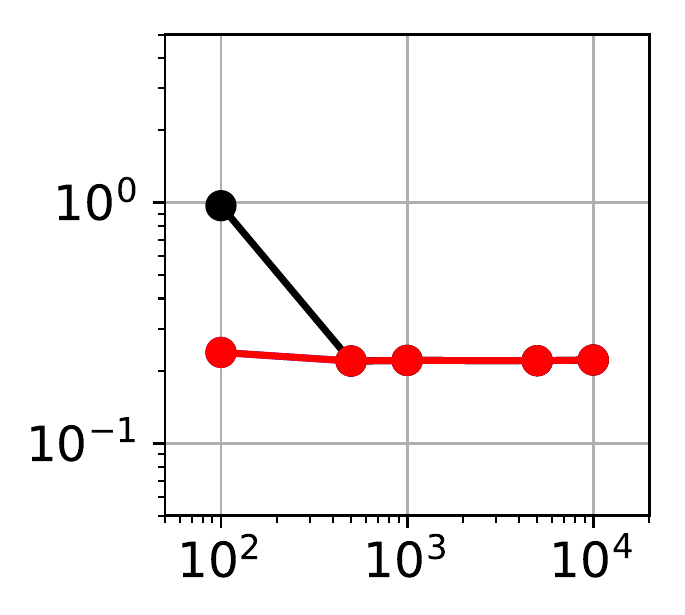}%
}
\subfloat[$n=1$, $N=4$]{
\includegraphics[width=0.24\textwidth]{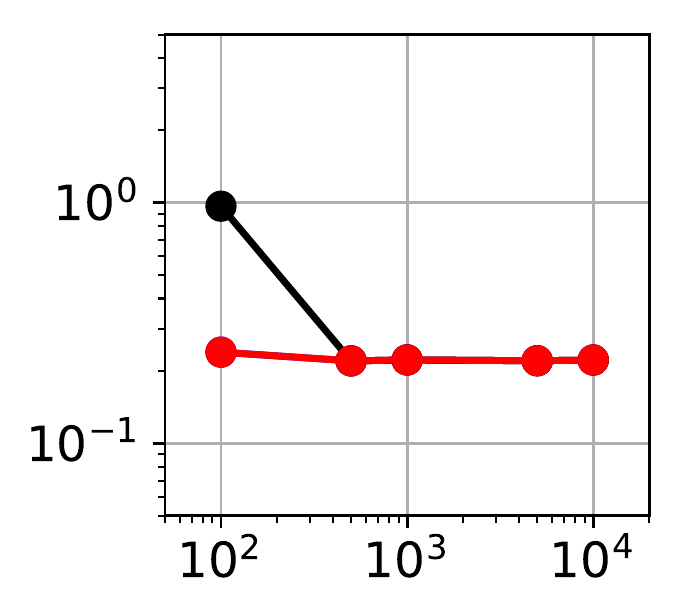}%
}
\subfloat[$n=1$, $N=5$]{
\includegraphics[width=0.24\textwidth]{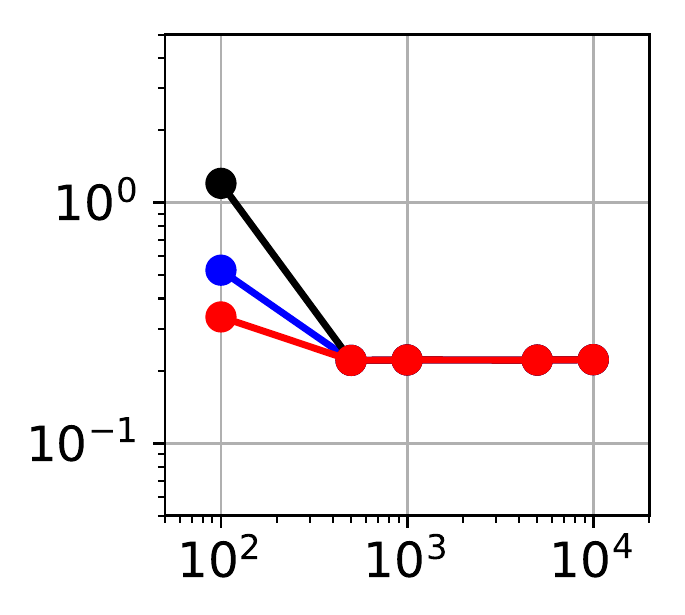}%
}
\\
\subfloat[$n=2$, $N=2$]{
\includegraphics[width=0.24\textwidth]{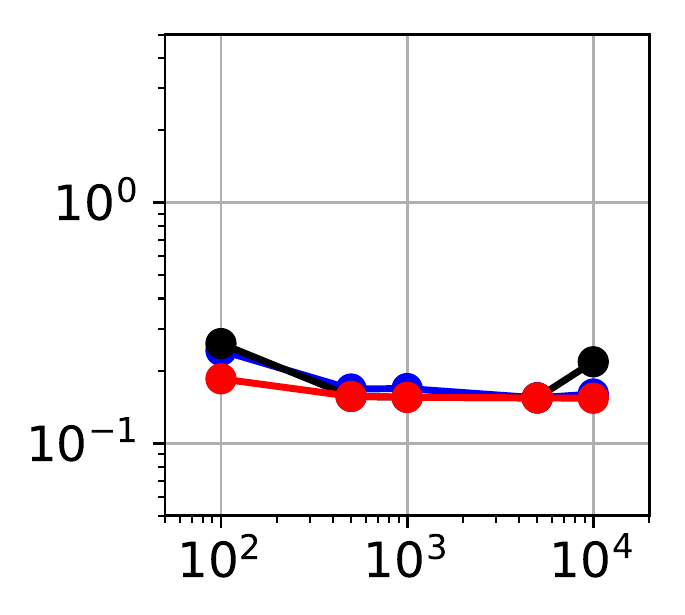}%
}
\subfloat[$n=2$, $N=3$]{
\includegraphics[width=0.24\textwidth]{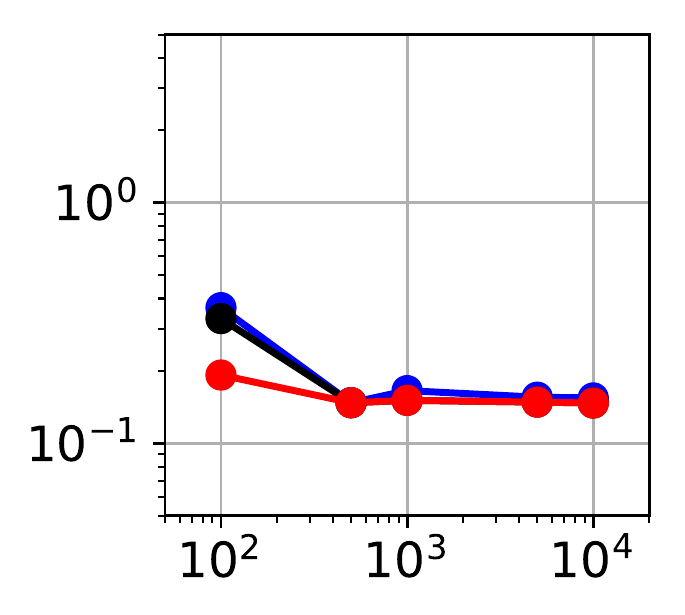}%
}
\subfloat[$n=2$, $N=4$]{
\includegraphics[width=0.24\textwidth]{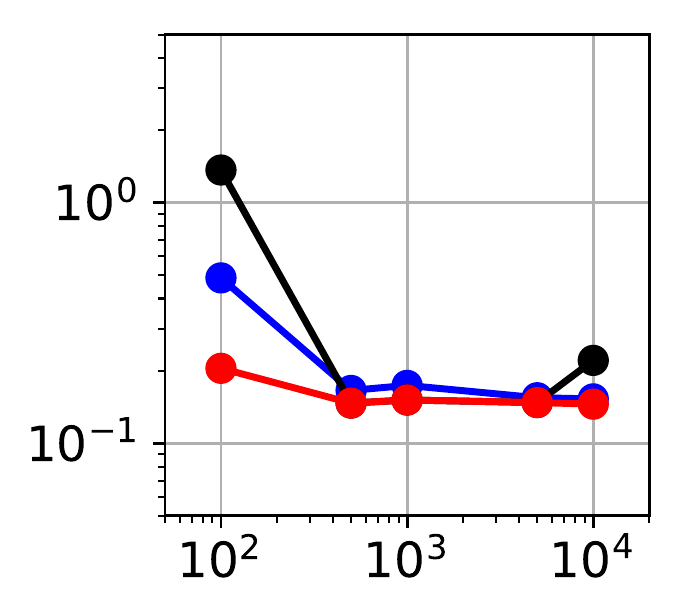}%
}
\subfloat[$n=2$, $N=5$]{
\includegraphics[width=0.24\textwidth]{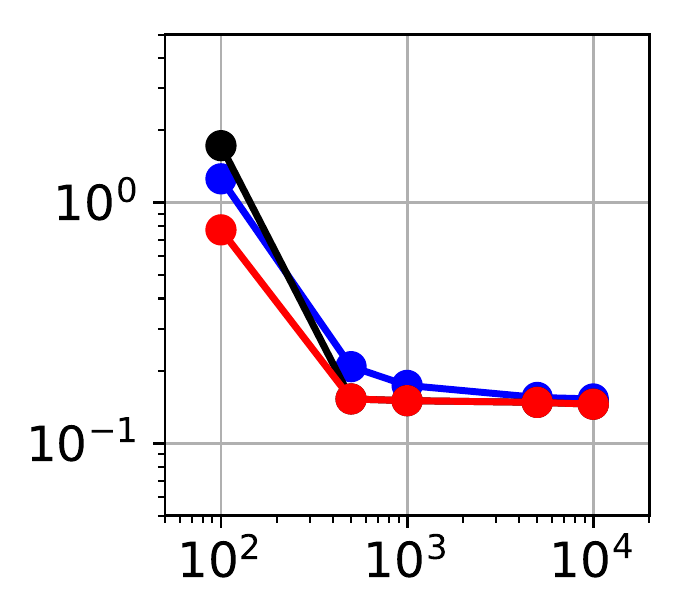}%
}
\\
\subfloat[$n=3$, $N=2$]{
\includegraphics[width=0.24\textwidth]{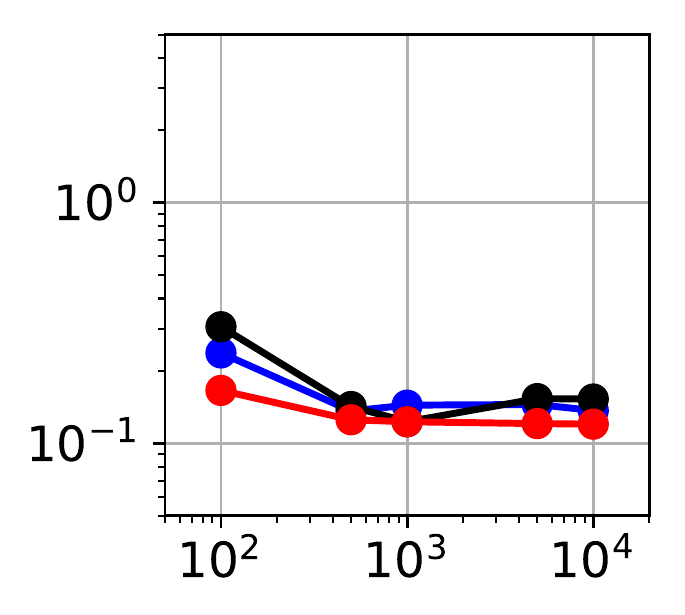}%
}
\subfloat[$n=3$, $N=3$]{
\includegraphics[width=0.24\textwidth]{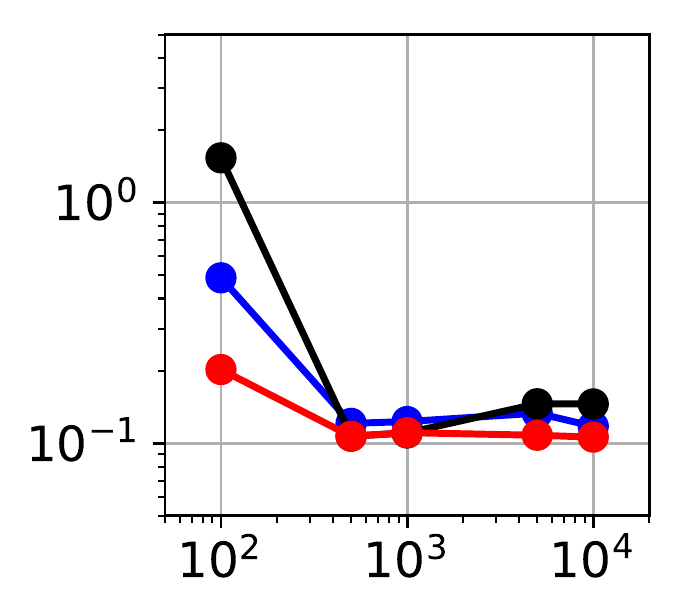}%
}
\subfloat[$n=3$, $N=4$]{
\includegraphics[width=0.24\textwidth]{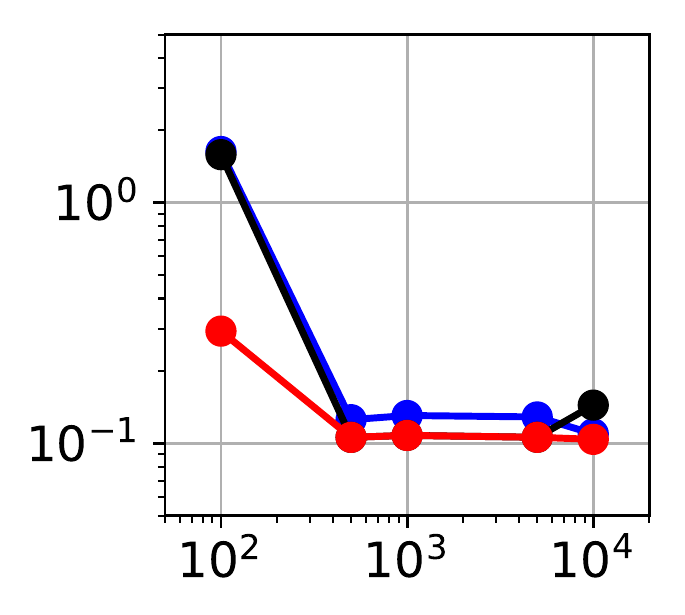}%
}
\subfloat[$n=3$, $N=5$]{
\includegraphics[width=0.24\textwidth]{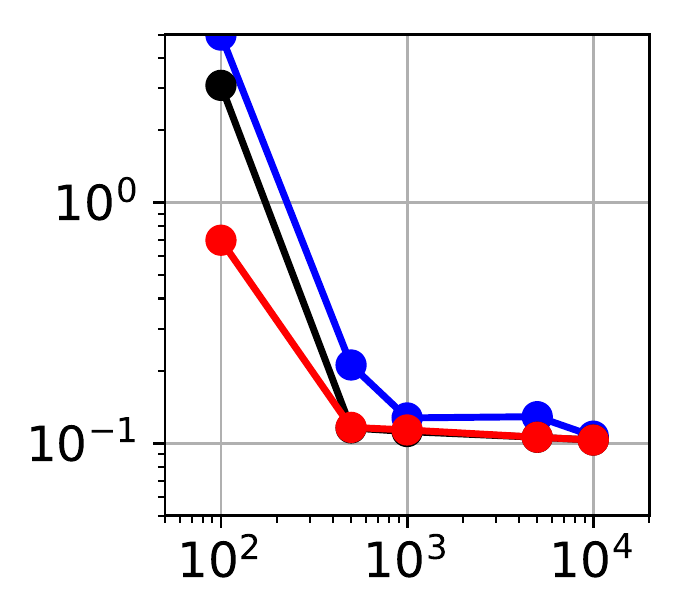}%
}
\\
\subfloat[$n=4$, $N=2$]{
\includegraphics[width=0.24\textwidth]{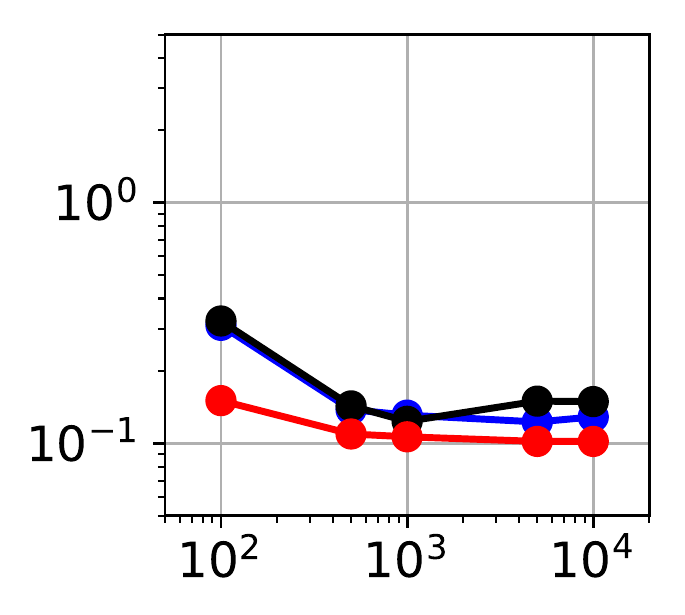}%
}
\subfloat[$n=4$, $N=3$]{
\includegraphics[width=0.24\textwidth]{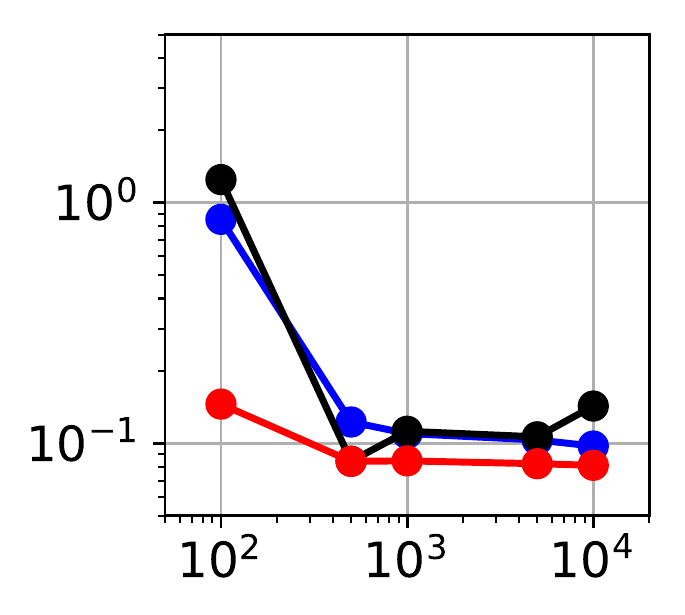}%
}
\subfloat[$n=4$, $N=4$]{
\includegraphics[width=0.24\textwidth]{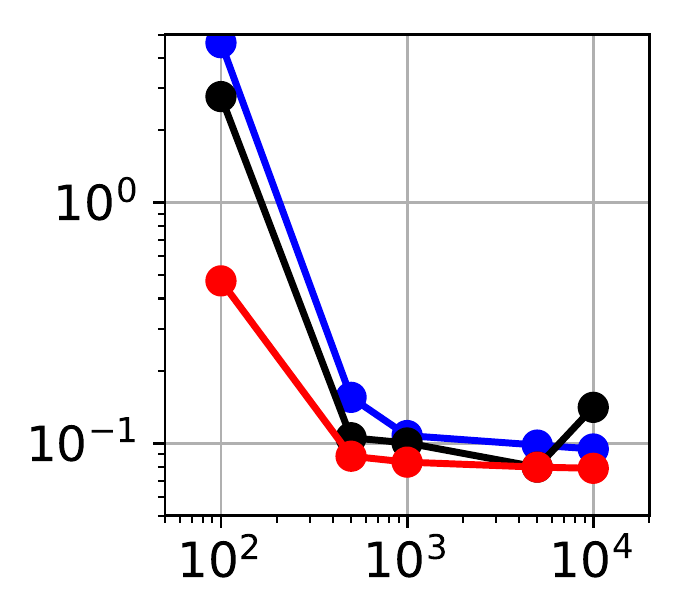}%
}
\subfloat[$n=4$, $N=5$]{
\includegraphics[width=0.24\textwidth]{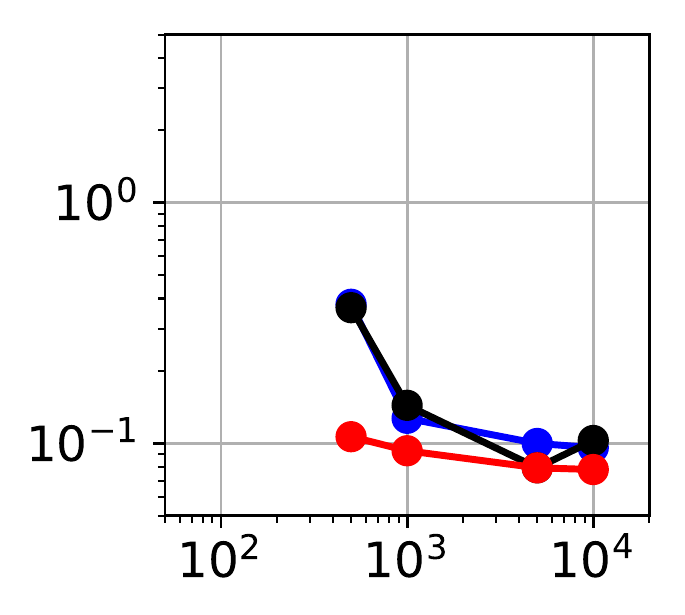}%
}
\caption{Average pointwise relative errors in the polynomial ridge approximation---fitted with 20 iterations of the alternating scheme in Algorithm \ref{alg:alternate}---on a testing set of 10000 samples as a function of the number of training samples.}
\label{fig:errs}
\end{figure}

For completeness, we study the error in the fitted ridge approximation as the number of training samples increases, $M\in\{100, 500, 1000, 5000, 10000\}$. We emphasize that it is not our primary goal to show that the ridge approximation model is ideal for the drag quantity of interest. To estimate error, we split the 20000-sample Latin hypercube design into a training set of maximum size 10000 and a testing set with 10000 samples. The testing error is computed as the average relative pointwise error over the testing set,
\begin{equation}
\frac{1}{M_{\text{test}}} \sum_{j=1}^{M_{\text{test}}} \frac{| f(\vx_j) - p_N(\mU^T\vx_j,\theta) |}{| f(\vx_j) |},
\end{equation}
where $\mU$ and $\theta$ are fitted with 20 iterations of the alternating heuristic from Algorithm \ref{alg:alternate}. Figure \ref{fig:errs} shows the testing error as a function of the number $M$ of training samples for each combination of $N$ (polynomial degree) and $n$ (number of linear combinations). Note that 100 samples is too few to fit a polynomial of degree $N=5$ in $n=4$ variables, so the bottom right subfigure is missing one data point relative to the other subfigures. The colors indicate the initial subspace guess in Algorithm \ref{alg:alternate}: black is the first $n$ columns of the identity matrix, blue is the average error over 10 random $n$-dimensional subspaces, and red is the active subspace. In most cases, the error is slightly smaller with the active subspace initial guess. However, all initial guesses lead to roughly the same testing error, which flattens after 500 training samples. Also, increasing the polynomial degree beyond $N=2$ has hardly any impact on the testing error. Increasing the number of linear combinations has a minor effect---decreasing the testing error by approximately 50\% from $n=1$ to $n=4$. 

\begin{figure}[!h]
\centering
\subfloat[Train]{
\includegraphics[width=0.48\textwidth]{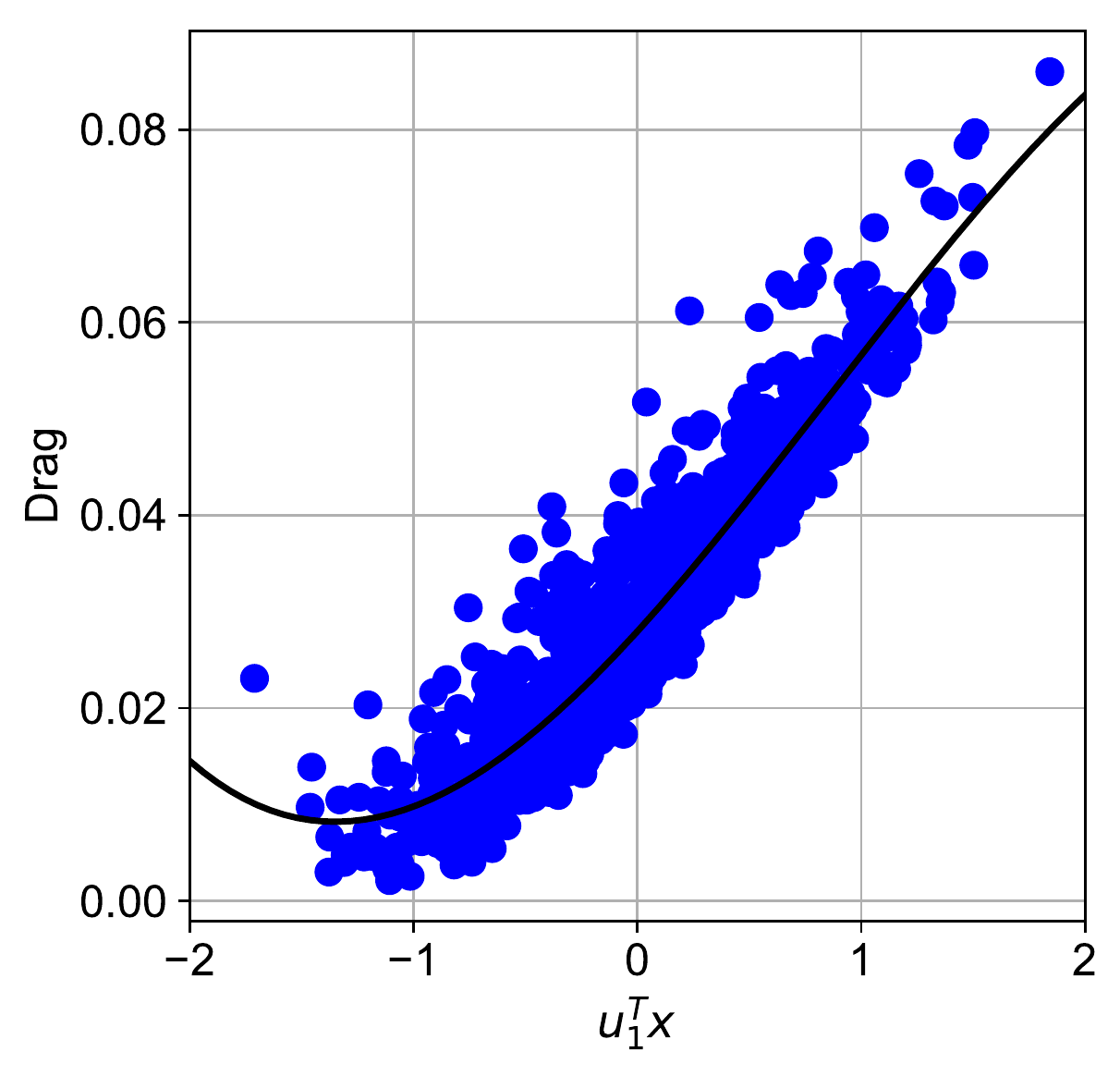}%
}
\subfloat[Test]{
\includegraphics[width=0.48\textwidth]{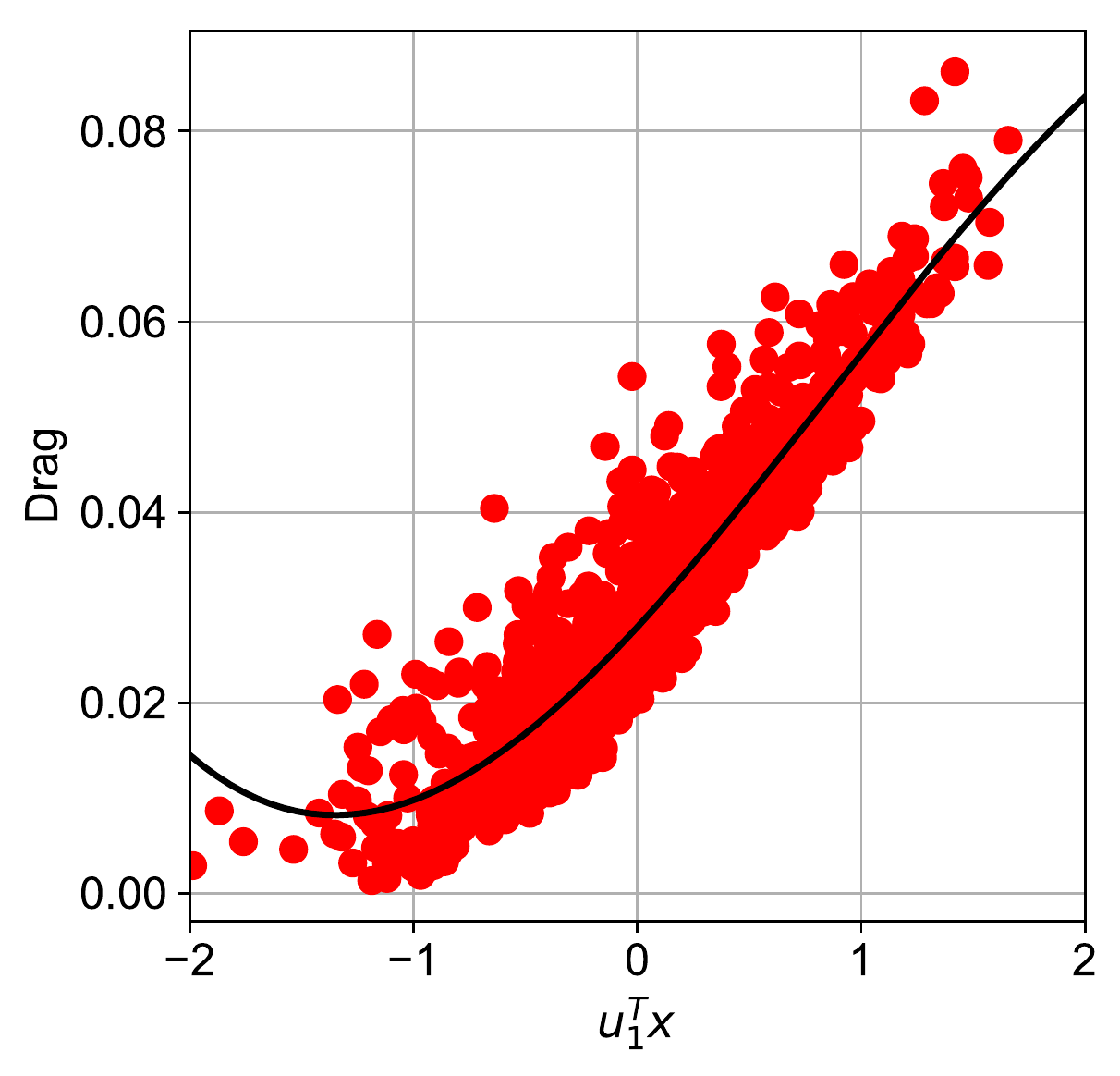}%
}
\caption{One-dimensional shadow plots on a training (left) and testing (right) set of 1000 input/output pairs. The black line is the fitted univariate polynomial of degree $N=5$ used in the ridge approximation. The globally monotonic structure may yield insight for specific design problems.}
\label{fig:1dshadow}
\end{figure}

\begin{figure}[!h]
\centering
\subfloat[Train]{
\includegraphics[width=0.48\textwidth]{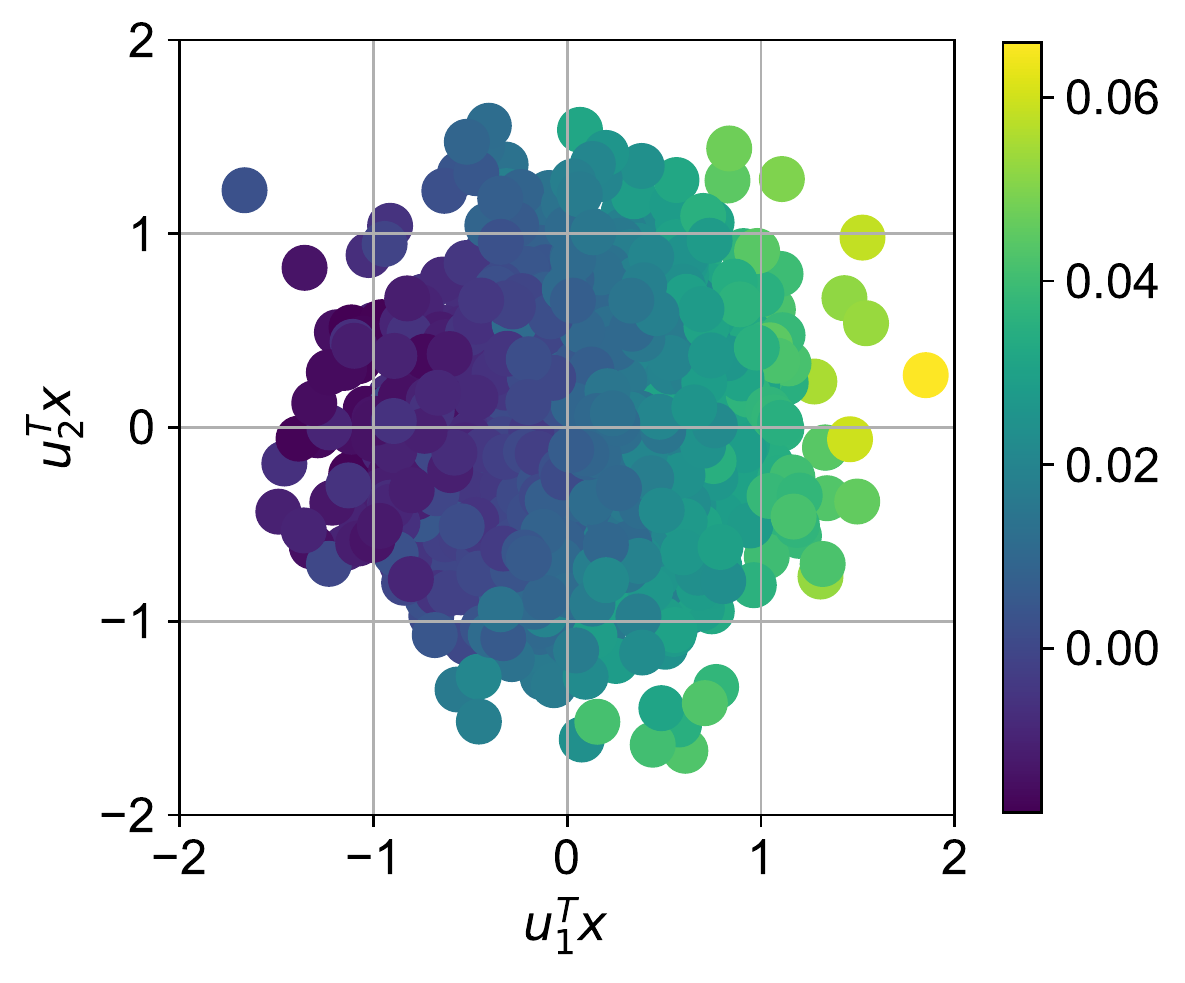}%
}
\subfloat[Test]{
\includegraphics[width=0.48\textwidth]{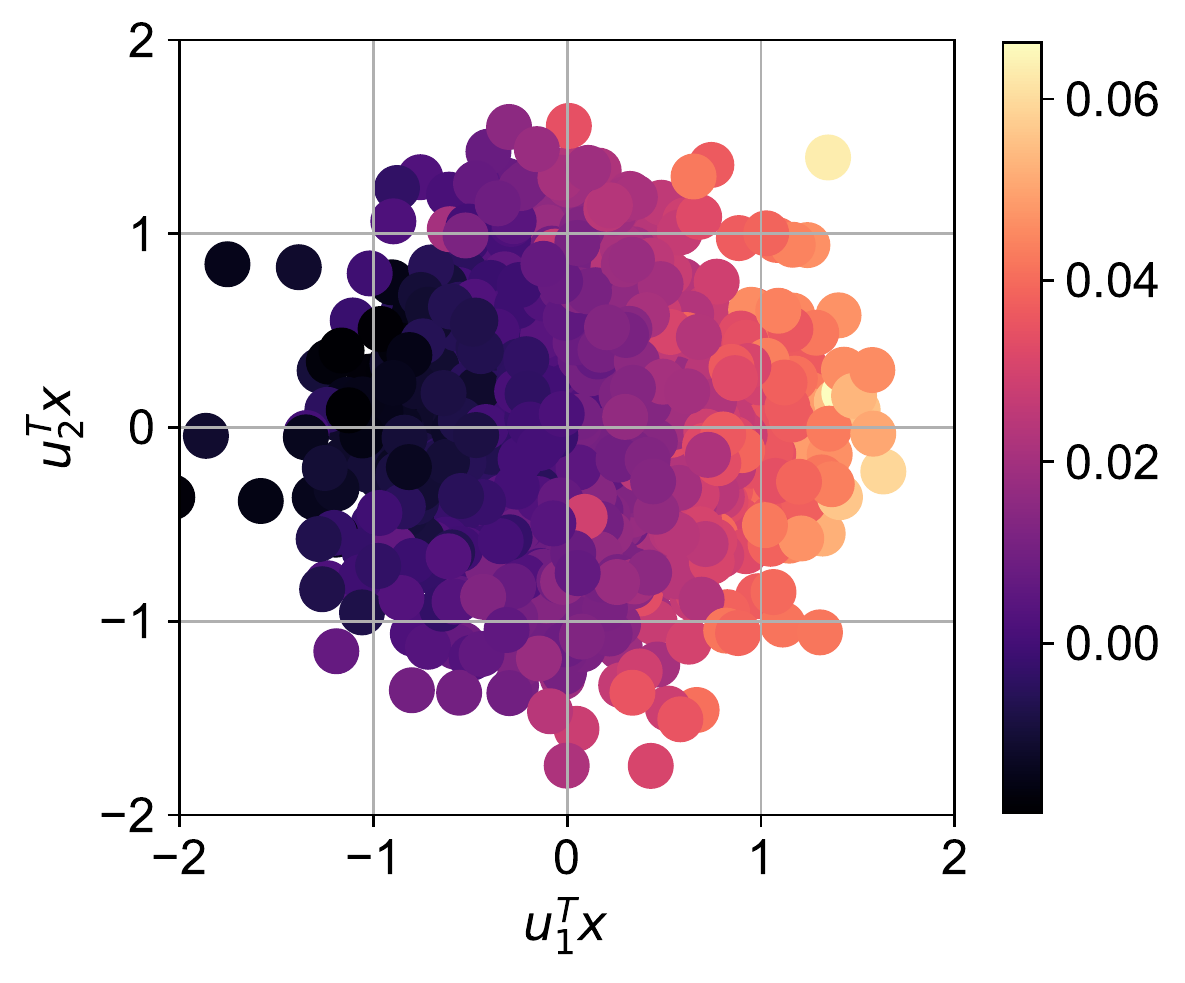}%
}
\\
\subfloat[Contour]{
\includegraphics[width=0.48\textwidth]{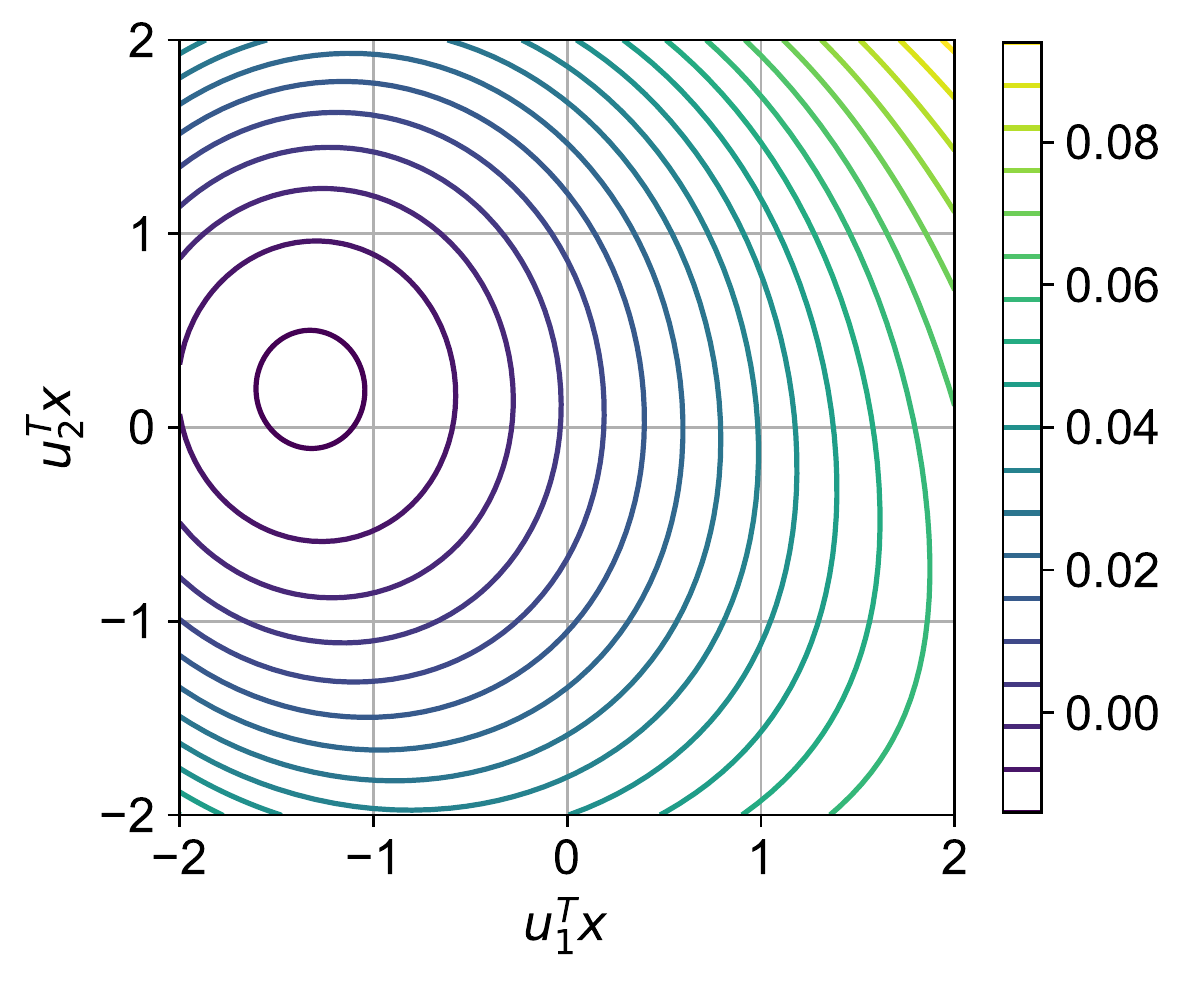}%
}
\caption{Two-dimensional shadow plots on a training (top left) and testing (top right) set of 1000 input/output pairs. The contour plot (bottom) is the bivariate polynomial of degree $N=5$ used in the ridge approximation.}
\label{fig:2dshadow}
\end{figure}

One advantage of the ridge approximation with $n=1$ or $n=2$ linear combinations is the possibility of visualizing the input/output relationship with simple scatter plots. More precisely, for the data set of input/output pairs $(\vx_i, f(\vx_i))$, if $\mU$ has 1 or 2 columns, then we can create scatter plots of the pairs $(\mU^T\vx_i, f(\vx_i))$. In regression problems where the data set consists of predictor/response pairs from an unknown joint density, such plots are called \emph{sufficient summary plots}~\cite{Cook1998}. The name \emph{sufficient} is used to imply that the plots are statistically sufficient for the data set, which is a technical condition on the linear combination weights $\mU$; for details, see Cook~\cite{Cook1998}. Since our data sets are not random in the sense of regression data (i.e., the function $f(\vx)$ is deterministic), we call these plots \emph{shadow plots}, since they resemble a shadow of a high-dimensional surface projected on a wall. For the $n=1$ case, in-to and out-of the page represents $m-1$-dimensional subspaces in the input space. 

Figure \ref{fig:1dshadow} shows the one-dimensional shadow plots of 1000 training samples and 1000 testing samples when $n=1$ along with the one-dimensional polynomial ridge approximation with degree $N=5$ (black line). The similarity between the plots suggests that the ridge direction $\mU$ would be similar if computed from the testing samples. The plot can be used to assess the quality of the one-dimensional ridge approximation. By visual inspection, the error in the one-dimensional ridge approximation may be up to roughly 25\% for this problem, and this is consist with the top row of plots in Figure \ref{fig:errs}. Although this error may be unacceptably large for pointwise approximation purposes (e.g., response surface-based optimization or uncertainty quantification), the plot reveals a globally monotonic trend in drag as a function of the shape parameters over the parameter space. Such insight may be exploited when designing an airfoil for drag. Moreover, the plot may assist in determining the set of shapes that produce sufficiently small drag for a constrained optimization; Constantine et al.~\cite{Constantine2015exploiting} use a similar approach for quantifying uncertainty in a scramjet model. Figure \ref{fig:2dshadow} shows the two-dimensional shadow plots, where the color is the value of drag. The contour plot shows the contours of the bivariate polynomial defining the two-dimensional ridge approximation. The quadratic behavior is apparent in the bivariate polynomial; again, this insight may be valuable to designers.

The results of this experiment provide numerical support that the first $n$ eigenvectors of $\mC$ from \eqref{eq:C}---even a cheap estimate---provide a good initial subspace $\mU_0$ for the alternating heuristic for ridge approximation of this aerospace model. This enables us to fit a ridge approximation with a relatively high polynomial degree along important directions in the model's input space.

%% file: sec4-conclusion.tex
\section{Summary and conclusions}

Motivated by response surface construction for expensive computational models with several input parameters, we study ridge approximation for functions of several variables. A ridge function is constant along a set of directions in its domain, and the approximation problem is to find (i) optimal directions and (ii) an optimal function of the linear combinations of variables. For a fixed set of directions, the best approximation in the mean-squared sense is a particular conditional average. We define an optimal subspace as one that minimizes a mean-squared cost function over the Grassmann manifold of subspaces. We then prove that a particular subspace---the active subspace defined by the function's gradient---is a near-stationary point for an optimization defining the optimal subspace. We offer a heuristic to exploit this fact when fitting a ridge approximation. Our first numerical example shows a simple case where this heuristic fails; this case reveals a type of function for which the heuristic is ill-suited, namely, functions that oscillate rapidly along one direction while varying slowly but consistently along another. Our second numerical example demonstrates this heuristic's success with a polynomial-based alternating scheme to fit a ridge approximation applied to an aerospace design model with 18 parameters. The alternating scheme with the active subspace as the initial guess outperforms the same scheme with random initial subspaces. Given the prevalence of anisotropic parameter dependence in most complex physical simulations, we expect that ridge functions are appropriate forms for response surfaces approximations. The analyses and heuristics we present advance the state-of-the-art in ridge approximations.

%% file: sec-appendix3.tex
\section{Proof of Theorem \ref{thm:necsuff}}
\label{sec:app1}

%\noindent Let $\mC$, $\mW_1$, $\mW_2$, $\Lambda_1$, and $\Lambda_2$ be defined as in Section \ref{sec:stationary}. Note that
%\begin{equation}
%\label{eq:p1}
%\begin{aligned}
%\int \left\|\mW_2^T\nabla f(\vx)\right\|^2\,\rho(\vx)\,d\vx
%&= 
%\int \nabla f(\vx)^T\mW_2\mW_2^T\nabla f(\vx)\,\rho(\vx)\,d\vx\\
%&= 
%\int \tr\left(\mW_2^T\nabla f(\vx) \nabla f(\vx)^T\mW_2\right) \,\rho(\vx)\,d\vx\\
%&= 
%\tr\left(\mW_2^T\left( \int \nabla f(\vx) \nabla f(\vx)^T\,\rho(\vx)\,d\vx\right)\mW_2\right) \\
%&= 
%\tr\left(\mW_2^T\mC\mW_2\right) \\
%&= 
%\tr\left(\Lambda_2\right) \\
%&= \lambda_{n+1} + \cdots + \lambda_m.
%\end{aligned}
%\end{equation}
%Assume that $f(\vx)$ is constant along directions corresponding to the columns of  $\mW_2$. Then the directional derivatives $\mW_2^T\nabla f(\vx)$ are zero for all $\vx$. Then by \eqref{eq:p1}, $\lambda_{n+1}=\cdots=\lambda_m=0$.

%Next, assume $\lambda_{n+1}=\cdots=\lambda_m=0$. Since the integrand in \eqref{eq:p1} is a norm and $f(\vx)$ is differentiable, $\mW_2^T\nabla f(\vx)$ is zero for all $\vx$. Therefore, $f(\vx)$ is constant along directions corresponding to the columns of $\mW_2$. 

\noindent Let $\vw\in\mathbb{R}^m$ be in the null space of $\mC$ from \eqref{eq:C}, i.e., $\mC\vw=\mathbf{0}$. Then 
\begin{equation}
0 
\;=\;
\vw^T\,\mC\,\vw
\;=\; 
\int \big(\vw^T\nabla f(\vx)\big)^2\,\rho(\vx)\,d\vx.
\end{equation}
The integrand is the squared directional derivative of $f$ along $\vw$. Since the squared quantity is non-negative, its average equalling zero---combined with continuity of $f$ from Assumption \ref{assump1}---implies that the directional derivative $\vw^T\nabla f(\vx)$ is zero for all $\vx$ in the support of $\rho$. Therefore, $f$ is constant along $\vw$. 

Now assume that $f$ is constant along $\vw$ in the support of $\rho$. Then $\vw^T\nabla f(\vx)=0$ for all such $\vx$. Then
\begin{equation}
0
\;=\; 
\int \big(\vw^T\nabla f(\vx)\big)^2\,\rho(\vx)\,d\vx
\;=\;
\vw^T\,\mC\,\vw,
\end{equation}
which implies that $\vw$ is in the null space of $\mC$, since $\mC$ is positive semidefinite.

\section{Proof of Theorem \ref{thm:gradbased}}
\label{sec:app2}
%Express $\mu$ as a function of $\vx$ and the complement subspace $\mV$,
%\begin{equation}
%\begin{aligned}
%\mu(\vx,\mV) 
%&= 
%\int f(\mU\mU^T\vx + \mV\vz)\,\pi(\vz)\,d\vz\\
%&= 
%\int f((\mI-\mV\mV^T)\vx + \mV\vz)\,\pi(\vz)\,d\vz.
%\end{aligned}
%\end{equation}
%Then write $R$ as a function of $\mV$,
%\begin{equation}
%R(\mV) \;=\; \frac{1}{2}\int (f(\vx)-\mu(\vx,\mV))^2\,\rho(\vx)\,d\vx.
%\end{equation}
\noindent For $R=R(\mV)$ from \eqref{eq:minV}, consider the gradient of $R$ on the Grassmann manifold $\mathbb{G}(m-n,n)$. 
\begin{equation}
\label{eq:gradR}
\begin{aligned}
\bnabla R(\mV) 
&= 
\bnabla \left(\frac{1}{2} \int (f(\vx)-\mu(\vx,\mV))^2\,\rho(\vx)\,d\vx\right)\\
&= 
\frac{1}{2} \int \bnabla (f(\vx)-\mu(\vx,\mV))^2\,\rho(\vx)\,d\vx\\
&= 
\int (f(\vx)-\mu(\vx,\mV))\,\bnabla (f(\vx)-\mu(\vx,\mV))\,\rho(\vx)\,d\vx\\
&= 
\int (f(\vx)-\mu(\vx,\mV))\,(\underbrace{\bnabla f(\vx)}_{=\,0} - \bnabla\mu(\vx,\mV))\,\rho(\vx)\,d\vx\\
&= 
\int (\mu(\vx,\mV)-f(\vx))\, \bnabla\mu(\vx,\mV)\,\rho(\vx)\,d\vx.
\end{aligned}
\end{equation}
Let $\mu_{ij}'$ be the $(i,j)$ element of $\bnabla \mu$, with $i=1,\dots,m$ and $j=1,\dots,m-n$. Using Cauchy-Schwarz, we can bound
\begin{equation}
\int (\mu-f)\,\mu_{ij}'\,\rho\,d\vx
\;\leq\;
\left(\int (\mu-f)^2\,\rho\,d\vx\right)^{\frac{1}{2}}
\left(\int (\mu_{ij}')^2\,\rho\,d\vx\right)^{\frac{1}{2}}.
\end{equation}
Then
\begin{equation}
\label{eq:frobound}
\begin{aligned}
\|\bnabla R(\mV)\|_F^2
&=
\sum_{i=1}^{m}\sum_{j=1}^{m-n}
\left(\int (\mu-f)\,\mu_{ij}'\,\rho\,d\vx\right)^2\\
&\leq
\sum_{i=1}^{m}\sum_{j=1}^{m-n} 
\left(\int (\mu-f)^2\,\rho\,d\vx\right)
\left(\int (\mu_{ij}')^2\,\rho\,d\vx\right)\\
&=
\left(\int (\mu-f)^2\,\rho\,d\vx\right)
\left(\int 
\sum_{i=1}^{m}\sum_{j=1}^{m-n} 
(\mu_{ij}')^2\,\rho\,d\vx\right)\\
&=
\left(\int (\mu-f)^2\,\rho\,d\vx\right)
\left(\int \|\bnabla \mu\|_F^2\,\rho\,d\vx\right).
\end{aligned}
\end{equation}
Recall Edelman's formula for the Grassmann gradient~\cite[Section 2.5.3]{Edelman1998},
\begin{equation}
\bnabla \mu(\vx,\mV) \;=\; (\mI-\mV\mV^T)\,\ppV \mu(\vx,\mV)
\;=\; \mU\mU^T\ppV \mu(\vx,\mV),
\end{equation}
where $\ppV\mu$ is the $m\times(m-n)$ matrix of partial derivatives of $\mu$ with respect to the elements of $\mV$. For Gaussian $\rho$, the conditional density
\begin{equation}
\pzgy \;=\; \pi(\vz) \;\propto\; \exp\left(\frac{-\vz^T\vz}{2}\right)
\end{equation}
is independent of $\mV$. Therefore,
\begin{equation}
\begin{aligned}
\ppV \mu(\vx,\mV)
&=
\ppV \int f((\mI-\mV\mV^T)\vx + \mV\vz)\,\pi(\vz)\,d\vz\\
&=
\int \ppV f((\mI-\mV\mV^T)\vx + \mV\vz)\,\pi(\vz)\,d\vz.
\end{aligned}
\end{equation}
Next we examine the gradient of $f$ with respect to the elements of $\mV$. For notation, define $\vs$ as
\begin{equation}
\vs \;=\;
\vs(\vx,\vz,\mV) \;=\;
(\mI-\mV\mV^T)\vx + \mV\vz.
\end{equation}
Let $v_{ij}$ be the $(i,j)$ element of $\mV$, and compute the derivative,
\begin{equation}
\ppvij f(\vs) \;=\;
\nabla f(\vs)^T\left(\ppvij \vs \right).
\end{equation}
The derivative of $\vs$ is 
\begin{equation}
\begin{aligned}
\ppvij \vs 
&=
\ppvij \left(
(\mI -\mV\mV^T)\vx + \mV\vz
\right)\\
&= \ve_i\,\vv_j^T\vx + x_i\,\vv_j + \ve_i\,z_j,
\end{aligned}
\end{equation}
where $\ve_i$ is the $i$th column of the $m\times m$ identity matrix, $\vv_j$ is the $j$th column of $\mV$, $x_i$ is the $i$th component of $\vx$, and $z_j$ is the $j$th component of $\vz$. Then
\begin{equation}
\begin{aligned}
\ppvij f(\vs) 
&= 
f_i(\vs)\,\vv_j^T\vx + x_i\,\nabla f(\vs)^T\vv_j + f_i(\vs)\,z_j\\
&= 
(f_i(\vs)\,\vx^T + x_i\,\nabla f(\vs)^T)\vv_j + f_i(\vs)\,z_j,
\end{aligned}
\end{equation}
where $f_i$ is the $i$th component of the gradient vector $\nabla f$. Putting $i$'s and $j$'s together,
\begin{equation}
\ppV f(\vs) 
\;=\;
\left(\nabla f(\vs)\,\vx^T + \vx\,\nabla f(\vs)^T\right)\,\mV + \nabla f(\vs)\,\vz^T.
\end{equation}
Then, for $f=f(\vs)$, $\pi=\pi(\vz)$, and $\nabla f=\nabla f(\vs)$,
\begin{equation}
\begin{aligned}
\int \ppV f\,\pi\,d\vz
&=
\int
\left(\nabla f\,\vx^T + \vx\,\nabla f^T\right)\,\mV + \nabla f\,\vz^T
\,\pi\,d\vz\\
&=
\left(\vg\,\vx^T + \vx\,\vg^T\right)\,\mV + \int\nabla f\,\vz^T\,\pi\,d\vz,
\end{aligned}
\end{equation}
where
\begin{equation}
\vg \;=\;
\vg(\vx) \;=\;
\int \nabla f((\mI-\mV\mV^T)\vx + \mV\vz)\,\pi(\vz)\,d\vz.
\end{equation}
By the Lipschitz continuity of $f$, $\|\vg\|\leq L$.
Then we can bound the norm of Grassmann gradient of $\mu$ as
\begin{equation}
\begin{aligned}
\|\bnabla \mu\|_F
&=
\left\| \mU\mU^T \ppV \mu \right\|_F\\
&\leq
\left\| \ppV \mu \right\|_F\\
&=
\left\| \int \ppV f\,\pi\,d\vz \right\|_F\\
&=
\left\| 
\left(\vg\,\vx^T + \vx\,\vg^T\right)\,\mV + \int\nabla f\,\vz^T\,\pi\,d\vz 
\right\|_F\\
&\leq
\left\| 
\left(\vg\,\vx^T + \vx\,\vg^T\right)\,\mV
\right\|_F + \left\| 
\int \nabla f\,\vz^T\,\pi\,d\vz
\right\|_F\\
&\leq
\left\| 
\vg\,\vx^T + \vx\,\vg^T
\right\|_F + \left(\int \| \nabla f\,\vz^T\|_F^2\,\pi\,d\vz\right)^{\frac{1}{2}}\\
&\leq
2\,\|\vg\|\,\|\vx\| + \left(\int \|\nabla f\|^2\|\vz\|^2\,\pi\,d\vz\right)^{\frac{1}{2}}\\
&\leq
2\,L\,\|\vx\| + \left(L^2\,\int \|\vz\|^2\,\pi\,d\vz\right)^{\frac{1}{2}}\\
&=
L\,\left(2\,\|\vx\| + (m-n)^{\frac{1}{2}}\right).
\end{aligned}
\end{equation}
Therefore,
\begin{equation}
\begin{aligned}
\int \|\bnabla \mu\|_F^2 \,\rho\,d\vx
&\leq
L^2 \int \left(2\,\|\vx\| + (m-n)^{\frac{1}{2}}\right)^2\,\rho\,d\vx \\
&\leq
L^2 \,\left(2m^{\frac{1}{2}} + (m-n)^{\frac{1}{2}}\right)^2.
\end{aligned}
\end{equation}
Combining this with \eqref{eq:frobound}, we have
\begin{equation}
\label{eq:0}
\|\bnabla R(\mV)\|_F 
\;\leq\;
L\,\left(2m^{\frac{1}{2}} + (m-n)^{\frac{1}{2}}\right)\,
\left(\int (\mu(\vx,\mV)-f(\vx))^2\,\rho(\vx)\,d\vx\right)^{\frac{1}{2}}.
\end{equation}
Note that the Poincar\'{e} constant for the Gaussian density is $C=1$~\cite{chen1982inequality}. Then combining \eqref{eq:0} with Theorem \ref{thm:approx} achieves the desired result.